\documentclass[a4paper]{article}
\usepackage{amsmath}\usepackage{epsf,amsfonts,amsthm}
\usepackage{fontenc,indentfirst, delarray,amsfonts,amsmath,amssymb}

\newcommand{\be}{\begin{equation}}
\newcommand{\ee}{\end{equation}}
\newcommand{\bea}{\begin{eqnarray*}}
\newcommand{\eea}{\end{eqnarray*}}
\newcommand{\w}{\wedge}
\newcommand{\p}{\partial}
\newcommand{\la}{\langle}
\newcommand{\ra}{\rangle}
\newcommand{\raa}{\rightarrow}
\newcommand{\Ci}{C^{\infty}}
\newcommand{\E}{\ell}
\newcommand{\N}{\mathbb{N}}

\newcommand{\R}{\mathbb{R}}

\newcommand{\C}{\mathbb{C}}
\newcommand{\Q}{\mathbb{Q}}

\newcommand{\lp}{\left(}
\newcommand{\rp}{\right)}

\newcommand{\op}[1]{\!\!\mathop{\rm ~#1}\nolimits}

\mathchardef\za="710B  %\alpha
\mathchardef\zb="710C  %\beta
\mathchardef\zg="710D  %\gamma
\mathchardef\zd="710E  %\delta
\mathchardef\zve="710F %\epsilon
\mathchardef\zz="7110  %\zeta
\mathchardef\zh="7111  %\eta
\mathchardef\zy="7112 %\theta
\mathchardef\zi="7113  %\iota
\mathchardef\zk="7114  %\kappa
\mathchardef\zl="7115  %\lambda
\mathchardef\zm="7116  %\mu
\mathchardef\zn="7117  %\nu
\mathchardef\zx="7118  %\xi
\mathchardef\zp="7119  %\pi
\mathchardef\zr="711A  %\rho
\mathchardef\zs="711B  %\sigma
\mathchardef\zt="711C  %\tau
\mathchardef\zu="711D  %\upsilon
\mathchardef\zf="711E %\phi
\mathchardef\zq="711F  %\chi
\mathchardef\zc="7120  %\psi
\mathchardef\zw="7121  %\omega
\mathchardef\ze="7122  %\varepsilon
\mathchardef\zvy="7123  %\vartheta
\mathchardef\zvw="7124  %\varomega
\mathchardef\zvr="7125 %\varrho
\mathchardef\zvs="7126 %\varsigma
\mathchardef\zvf="7127  %\varphi
\mathchardef\zG="7000  %\Gamma
\mathchardef\zD="7001  %\Delta
\mathchardef\zY="7002  %\Theta
\mathchardef\zL="7003  %\Lambda
\mathchardef\zX="7004  %\Xi
\mathchardef\zP="7005  %\Pi
\mathchardef\zS="7006  %\Sigma
\mathchardef\zU="7007  %\Upsilon
\mathchardef\zF="7008  %\Phi
\mathchardef\zW="700A  %\Omega

\newtheorem{theo}{Theorem}
\newtheorem{prop}{Proposition}

\newtheorem{defi}{Definition}

\begin{document}

\title{Dequantized Differential Operators between Tensor Densities\\ as Modules over the Lie Algebra of Contact Vector Fields}
\author{Ya\"el Fr\'egier$\mbox{ }^*$, Pierre Mathonet$\mbox{ }^\dag$, Norbert Poncin \footnote{University of Luxembourg, Campus Limpertsberg, Institute
of Mathematics, 162A, avenue de la Fa\"iencerie, L-1511 Luxembourg
City, Grand-Duchy of Luxembourg, E-mail: yael.fregier@uni.lu,
norbert.poncin@uni.lu. The research of N. Poncin was supported by
grant R1F105L10. This author also thanks the Erwin Schr\"odinger
Institute in Vienna for hospitality and support during his visits
in 2006 and 2007. $^{\dag}$University of Li\`ege, Sart Tilman,
Institute of Mathematics, B37, B-4000 Li\`ege, Belgium, E-mail:
P.Mathonet@ulg.ac.be.}}\maketitle

\begin{abstract}\noindent In recent years, algebras and modules of differential operators have been
extensively studied. Equivariant quantization and dequantization
establish a tight link between invariant operators connecting
modules of differential operators on tensor densities, and module
morphisms that connect the corresponding dequantized spaces. In
this paper, we investigate dequantized differential operators as
modules over a Lie subalgebra of vector fields that preserve an
additional structure. More precisely, we take an interest in
invariant operators between dequantized spaces, viewed as modules
over the Lie subalgebra of infinitesimal contact or projective
contact transformations. The principal symbols of these invariant
operators are invariant tensor fields. We first provide full
description of the algebras of such affine-contact- and
contact-invariant tensor fields. These characterizations allow
showing that the algebra of projective-contact-invariant operators
between dequantized spaces implemented by the same density weight,
is generated by the vertical cotangent lift of the contact form
and a generalized contact Hamiltonian. As an application, we prove
a second key-result, which asserts that the Casimir operator of
the Lie algebra of infinitesimal projective contact
transformations, is diagonal. Eventually, this upshot entails that
invariant operators between spaces induced by different density
weights, are made up by a small number of building bricks that
force the parameters of the source and target spaces to verify
Diophantine-type equations.

\end{abstract}

\\


\section{Introduction}

Equivariant quantization, in the sense of C. Duval, P. Lecomte,
and V. Ovsienko, developed as from 1996, see \cite{LMT},
\cite{LO}, \cite{DLO}, \cite{PL}, \cite{BM}, \cite{DO},
\cite{BHMP}, \cite{BM2}. This procedure requires equivariance of
the quantization map with respect to the action of a
finite-dimensional Lie subgroup $G\subset\op{Diff}(\R^n)$ of the
symmetry group $\op{Diff}(\R^n)$ of configuration space $\R^n$,
or, on the infinitesimal level, with respect to the action of a
Lie subalgebra of the Lie algebra of vector fields. Such
quantization maps are well-defined globally on manifolds endowed
with a flat $G$-structure and lead to invariant star-products,
\cite{LO},\cite{DLO}. Equivariant quantization has first been
studied on vector spaces, mainly for the projective and conformal
groups, then extended in 2001 to arbitrary manifolds, see
\cite{PL2}. In this setting, equivariance with respect to all
arguments and for the action of the group of all (local)
diffeomorphisms of the manifold (i.e. naturality in the sense of
I. Kol\'a\v{r}, P. W. Michor, and J. Slov\'ak, \cite{KMS}) has
been ensured via quantization maps that depend on (the projective
class of) a connection. Existence of such natural and projectively
invariant quantizations has been investigated in several works,
\cite{MB}, \cite{MR}, \cite{SH}.

From the very beginning, equivariant quantization and symbol
calculus, and classification issues in Representation Theory of
Algebras appeared as dovetailing topics, see \cite{LMT},
\cite{LO}, \cite{PM2}, \cite{BHMP}, \cite{NP}. In these works
differential operators between sections of vector bundles have
been studied and classified as modules over the Lie algebra of
vector fields. Except for \cite{PM}, the case of differential
operators as representations of a subalgebra of vector fields that
preserve some additional structure, was largely uninvestigated.
The origin of this paper is the classification problem of
differential operators on a contact manifold between tensor
densities of possibly different weights (in the frame of
equivariant quantization it is natural to consider linear
differential operators between densities rather than between
functions, as [even mathematical] quantization maps should be
valued in a space of operators acting on a Hilbert or preHilbert
space), as modules over the Lie subalgebra of contact vector
fields.

Let us give a rough description of our approach to the preceding
multilayer problem. Further details can be found below.
Projectively equivariant quantization establishes a tight
connection between the ``quantum level''-- classification of
differential operators as representations of the algebra of
contact vector fields, and the ``classical level''-- quest for
intertwining operators between the corresponding modules of
symbols over the subalgebra of infinitesimal projective contact
transformations. These morphisms (in the category of modules) have
(locally) again symbols and these are tensor fields. The principal
symbol map intertwines the natural actions on morphisms and tensor
fields. Hence, the principal symbol of any ``classical''
intertwining operator is an invariant tensor field. These
invariant fields can be computed. However, it turns out that the
obvious technique that should allow lifting invariant tensor
fields to ``classical'' module morphisms is not sufficient for our
purpose. The Casimir operator (of the representation of
infinitesimal projective contact transformations on symbols)
proves to be an efficient additional tool. Calculation of the
Casimir itself requires a noncanonical splitting of the module of
symbols. This decomposition has been elaborated in a separate
paper, see \cite{FMP}.

In the present work, we investigate the ``classical level''
problem, i.e. we study ``dequantized'' differential operators
between tensor densities as modules over infinitesimal contact
transformations.

The paper is self-contained and organized as follows.

In Section \ref{ContGeo}, we recall essential facts in Contact
Geometry, which are relevant to subsequent sections. We place
emphasis on global formul\ae, as till very recently most of the
results were of local nature.

Section \ref{InfProjContTrans} provides the whole picture related
with infinitesimal projective contact transformations. A good
understanding of these upshots is crucial, particularly as regards
the calculation, in Section \ref{InvTensFields}, of invariant
tensor fields, and in consideration of the computation of the
aforementioned Casimir operator, see Section \ref{CasimirSection}.

Coordinate-free approaches to differential operators, their
symbols, and all involved actions are detailed in Section
\ref{DiffOpSymbEct}. This material is of importance with respect
to the geometric meaning of several invariant tensor fields
constructed later.

In Section \ref{InvTensFields}, we give a full description of the
algebra of affine-contact-invariant tensor fields (local
investigation), see Theorem \ref{invaff}, and of the algebra of
contact-invariant tensor fields (global result), see Theorem
\ref{invcont}.

A third main upshot, based on the preceding Section, is the
assertion that the algebra of projective-contact-invariant
operators between symbol modules ``implemented by the same density
weight'', is generated by two basic operators, the vertical
cotangent lift of the contact form and a generalized contact
Hamiltonian, both introduced in \cite{FMP}, see Theorem
\ref{ClassSpInvSame}, Section \ref{InvOpSameDens}.

As an application of the aforenoted noncanonical splitting of the
module of symbols into submodules, see \cite{FMP}, of Section
\ref{InfProjContTrans}, and of Section \ref{InvOpSameDens}, we
prove in Section \ref{CasimirSection}, that the Casimir operator
of the canonical representation of the Lie subalgebra of
infinitesimal projective contact transformations on the mentioned
symbol space, with respect to the Killing form, is diagonal, see
Theorem \ref{CasDiagForm}.

Eventually, the computation of this Casimir operator---actually a
challenge by itself---allows showing that the quest for
projective-contact-invariant operators between symbol modules
``implemented by different density weights'', can be put down to
the search of a small number of invariant building blocks between
(smaller) eigenspaces, see Section \ref{InvOpDiffDens}. Further,
each such brick forces the parameters of the source and target
symbol modules to verify a Diophantine-type equation.

\section{Remarks on Contact Geometry}\label{ContGeo}

A contact structure on a manifold $M$ is a co-dimension 1 smooth
distribution $\xi$ that is completely nonintegrable. Such a
distribution is locally given by the kernel of a nowhere vanishing
$1$-form $\za$ defined up to multiplication by a never vanishing
function. Since $(d\za)(X,Y)=X(\za(Y))-Y(\za(X))-\za([X,Y]),$
where notations are self-explaining, {\it integrability} of $\xi$,
i.e. closeness of sections of $\xi$ under the Lie bracket of
vector fields, would require that $d\za$ vanish on vectors in
$\zx$. By complete {\it nonintegrability} we mean that $d\za$ is
nondegenerate in $\xi$, for any locally defining $1$-form $\za$.
It follows that contact manifolds are necessarily odd-dimensional.
Eventually we get the following definition.
\begin{defi}
A contact manifold $M$ is a manifold of odd dimension $2 n + 1$,
together with a smooth distribution $\xi$ of hyperplanes in the
tangent bundle of $M$, such that $\alpha\wedge (d\alpha)^n$ is a
nevervanishing top form for any locally defining $1$-form $\za$.
Distribution $\xi$ is a contact distribution or contact structure
on $M$. A Pfaffian manifold $($A. Lichnerowicz$)$ or coorientable
contact manifold is a manifold $M$ with odd dimension $2n+1$
endowed with a globally defined differential $1$-form $\za$, such
that $\za\w(d\za)^n$ is a volume form of $M$. Form $\za$ $($which
defines of course a contact distribution on $M$$)$ is called a
contact form on $M$.
\end{defi}

\noindent{\bf Example 1}. Let
$(p_1,\ldots,p_n,q^1,\ldots,q^n,t,\tau)$ be canonical coordinates
in $\R^{2n+2}$, let $i:\R^{2n+1}\hookrightarrow\R^{2n+2}$ be the
embedding that identifies $\R^{2n+1}$ with the hyperplane $\tau=1$
of $\R^{2n+2}$, and let $\zs$ be the Liouville $1$-form of
$\R^{2n+2}$ (which induces the canonical symplectic structure of
$\R^{2n+2}$). It is easily checked that the pullback
$$\za=i^*\zs=\frac{1}{2}(\sum_{k=1}^n(p_k
dq^k - q^k dp_k)-dt)$$ of the Liouville form $\zs$ by embedding
$i$ is a contact form on $\R^{2n+1}$. Any coorientable contact
manifold $(M,\za)$ can locally be identified with $(\R^{2 n +
1},i^*\zs)$, i.e. Darboux' theorem holds true for contact manifolds.\\

\noindent{\bf Remark 1}. The preceding extraction of a contact
structure from a symplectic structure is the shadow of a tight
connection between contact and symplectic manifolds. If $(M,\za)$
is a coorientable contact manifold, if $\zp:M\times\R\raa M$ is
the canonical projection, and $s$ a coordinate function in $\R$,
the form $\zw=d(e^{s}\zp^*\za)$ is a symplectic form on
$M\times\R$, which is homogeneous with respect to $\p_{s}$, i.e.
$L_{\p_{s}}\zw=\zw$. This symplectic homogeneous manifold
$(M\times\R,\zw,\p_{s})$ is known as the {\it symplectization} the
initial contact manifold $(M,\za)$. Actually, there is a
$1$-to-$1$ correspondence between coorientable contact structures
on $M$ and homogeneous symplectic structures on $M\times\R$ (with
vector field $\p_{s}$). This relationship extends from the
contact-symplectic to the Jacobi-Poisson setting, see \cite{AL},
for super-Poissonization, see \cite{GU}. A coordinate-free
description of symplectization is possible. Consider a contact
manifold $(M,\xi)$, let $L\subset T^*M$ be the line subbundle of
the cotangent bundle, made up by all covectors that vanish on
$\xi$, and denote by $L_0$ the submanifold of $L$ obtained by
removing the $0$-section. The restriction to $L_0$ of the standard
symplectic form of $T^*M$ endows $L_0$ with a symplectic
structure, see \cite{VA, VO}. Eventually, a contact structure on a
manifold $M$ can be viewed as a line subbundle $L$ of the
cotangent bundle $T^*M$ such that the restriction to $L_0$ of the
standard symplectic form on $T^*M$ is symplectic. Of course, the
contact structure is coorientable if and only if $L$ is trivial.\\

\noindent{\bf Example 2}. Let $\za=i^*\zs$ be the standard contact
form of $\R^{2n+1}$. Set $x=(p,q,t,\tau)=(x',\tau)$. The open half
space $\R^{2n+2}_+=\{(x',\tau):\tau >0\}$, endowed with its
canonical symplectic structure $\zw$ and the Liouville vector
field $\zD=\frac{1}{2}{\cal E}$, where ${\cal E}$ is the usual
Euler field, can be viewed as symplectization of
$(\R^{2n+1},\za)$. Indeed, $\R^{2n+2}_+\simeq
\cup_{x'\in\R^{2n+1}} \{x'\}\times\{\tau(x',1):\tau>0\}$ is a line
bundle over $\R^{2n+1}$ with fiber coordinate $s=\op{ln}\tau^2$.
The projection of this bundle reads $\zp:\R^{2n+2}_+\ni
(x',\tau)\raa \tau^{-1}x'\in\R^{2n+1}.$ Clearly, $\zw$ has degree
$1$ with respect to $\zD$ and it is
easily checked that $\zw=d(e^s\zp^*\za)$ and $\zD=\p_s$.\\ %$p_s$ is tangent
% to the fiber at $s$, so it is $d_s(\sqrt{e^s}(x',1))$. When writing this
% vector not at $s=\tau(x',1)$, but at $(x',\tau)$, we find $\zD$.\\

In the following, unless otherwise stated, we consider
coorientable contact manifolds (or trivial line bundles).

\begin{defi} Let $M$ be a (coorientable) contact manifold. A
contact vector field is a vector field $X$ of $M$ that preserves
the contact distribution. In other words, for any fixed contact
form $\za$, there is a function $f_X\in\Ci(M)$, such that
$L_X\za=f_X\za$. We denote by $\op{CVect}(M)$ the space of contact
vector fields of $M$.
\end{defi}

It is easily seen that space $\op{CVect(M)}$ is a Lie subalgebra
of the Lie algebra $\op{Vect}(M)$ of all vector fields of $M$, but
not a $\Ci(M)$-module.\\

Let us now fix a contact form $\za$ on $M$ and view $d\za$ as a
bundle map $d\za:TM\raa T^*M$. It follows from the nondegeneracy
condition that the kernel $\op{ker}d\za$ is a line bundle and that
the tangent bundle of $M$ is canonically split:
$TM=\op{ker}\za\oplus\op{ker}d\za.$ Moreover,
\be\op{Vect(M)}=\op{ker}\za\oplus\op{ker}d\za,\label{CanSplit}\ee
where $\za$ and $d\za$ are now viewed as maps between sections. It
is clear that there is a unique vector field $E$, such that
$i_Ed\za=0$ and $i_E\za=1$ (normalization condition). This field
is called the {\it Reeb vector field}. It is {\it strongly
contact} in the sense that $L_E\za=0.$\\

Pfaffian structures, just as symplectic structures, can be
described by means of contravariant tensor fields. These fields
are obtained from $\za$ and $d\za$ via the musical map
$\flat:\op{Vect}(M)\ni X\raa (i_X\za)\za+i_Xd\za\in\zW^1(M)$,
which is a $\Ci(M)$-module isomorphism, see e.g. \cite{LLMP}. The
contravariant objects in question are the Reeb vector field
$E=\flat^{-1}(\za)\in\op{Vect}(M)=:{\cal X}^1(M)$ and the bivector
field $\zL\in{\cal X}^2(M)$, defined by
$$\zL(\zb,\zg)=(d\za)(\flat^{-1}(\zb),\flat^{-1}(\zg)),$$
$\zb,\zg\in\zW^1(M)$. They verify $[\zL,\zL]_{\op{SCH}}=2E\w\zL$
and $L_E\zL=0$, where $[.,.]_{\op{SCH}}$ is the Schouten-Nijenhuis
bracket. Hence, any (coorientable) contact manifold is a {\it
Jacobi manifold}.

Let us recall that Jacobi manifolds are precisely manifolds $M$
endowed with a vector field $E$ and a bivector field $\zL$ that
verify the two preceding conditions. The space of functions of a
Jacobi manifold $(M,\zL,E)$ carries a Lie algebra structure,
defined by \be\{h,g\}=\zL(dh,dg)+hEg-gEh,\label{JacBrack}\ee
$h,g\in\Ci(M)$. The Jacobi identity for this bracket is equivalent
with the two conditions $[\zL,\zL]_{\op{SCH}}=2E\w\zL$ and
$L_E\zL=0$ for Jacobi manifolds (these conditions can also be
expressed in terms of the Nijenhuis-Richardson bracket, see
\cite{NR}). It is well-known that the ``Hamiltonian map'' \be
X:\Ci(M)\ni h\raa
X_h=i_{dh}\zL+hE\in\op{Vect}(M)\label{JacHamMap}\ee is a Lie
algebra homomorphism: $X_{\{h,g\}}=[X_h,X_g]$. If $\op{dim}M=2n+1$
and $E\w\zL^n$ is a nowhere vanishing tensor field, manifold $M$
is coorientably contact. Furthermore, if we fix, in the Pfaffian
case, a contact form $\za$, we get a Lie algebra isomorphism \be
X:\Ci(M)\ni h\raa X_h\in\op{CVect}(M)\label{ContHamMap}\ee between
functions and contact vector fields, see \cite{VA}. It follows
from the above formul\ae$\mbox{}$ that $\za(X_h)=h$.

The main observation is that Jacobi brackets, see
(\ref{JacBrack}), are first order bidifferential operators. This
fact is basic in many recent papers, see e.g. \cite{GM} (inter
alia for an elegant approach to graded Jacobi cohomology), or
\cite{GU} (for Poisson-Jacobi reduction).\\ % Note also that in
% \cite{JG}, J. Grabowski suggests investigating Jacobi structures
% for sections of arbitrary line bundles rather than for functions,
% i.e. sections of a trivial line bundle. The trivial case has been
% studied by A. A. Kirillov, \cite{AK}, who proved
% that any local Lie structure on sections of a trivial line bundle is a Jacobi structure on the base.\\

After the above global formul\ae$\mbox{}$ and fundamental facts on
Contact Geometry, we continue with other remarks that are of
importance for our investigations. The setting is still a
($2n+1$)-dimensional contact manifold $M$ with fixed contact form
$\za$. Contraction of the equation $L_{X_h}\za=f_{X_h}\za$,
$h\in\Ci(M)$, with the Reeb field $E$ leads to $f_{X_h}=E(h)$. If
$\zW$ denotes the volume $\zW=\za\w(d\za)^n$, it is clear that,
for any contact vector field $X$, we have $L_X\zW=(n+1)f_X\;\zW$.
Hence, \be\op{div}_{\zW}X=(n+1)f_X, \forall
X\in\op{CVect}(M),\label{DivContact}\ee and
$\op{div}_{\zW}X_h=(n+1)E(h)$, for any $h\in\Ci(M).$ It follows
that for all $h,g\in\Ci(M)$,
\be\{h,g\}=X_h(g)-gE(h)=X_h(g)-\frac{1}{n+1}g\op{div}_{\zW}X_h=L_{X_h}\tilde{g},\label{LieDensities}\ee
where $\tilde{g}$ is function $g$ viewed as tensor density of
weight $-1/(n+1)$. Tensor densities will be essential below. For
details on densities, we refer the reader to \cite{FMP}. The
afore-depicted Lie algebra isomorphism $X$ between functions and
contact vector fields, is also a $\op{CVect}(M)$-module
isomorphism, if we substitute the space ${\cal
F}_{\frac{-1}{n+1}}(M)$ of tensor densities of weight $-1/(n+1)$
for the space of functions (of course, the contact action is
$L_{X_{\tilde h}}{\tilde g}=\{{\tilde h},{\tilde g}\}$ on
densities, and it is the adjoint action on contact fields). Note
that this distinction between functions and densities is necessary
only if the module structure is concerned.\\

We now come back to splitting (\ref{CanSplit}). If we denote by
$\op{TVect(M)}$ the space of tangent vector fields, i.e. the space
$\op{ker}\za$ of those vector fields of $M$ that are tangent to
the contact distribution, this decomposition also reads
$$\op{Vect}(M)=\op{TVect(M)}\oplus\;\Ci(M)E.$$ As abovementioned,
our final goal is the solution of the multilayer classification
problem of differential operators between tensor densities on a
contact manifold, as modules over the Lie algebra of contact
vector fields. This question naturally leads to the quest for a
splitting of some $\op{CVect}(M)$-modules or
$\op{sp_{(2n+2)}}$-modules of symbols, see below, and in
particular of the module $\op{Vect(M)}$ itself. Space
$\op{TVect}(M)$, which is of course not a Lie algebra, is a
$\Ci(M)$-module and a $\op{CVect}(M)$-module. The last upshot
follows directly from formula $i_{[X,Y]}=[L_X,i_Y]$,
$X,Y\in\op{Vect(M)}$. The second factor $\Ci(M)E$ however, is
visibly not a $\op{CVect(M)}$-module (for instance
$[E,X_h]=[X_1,X_h]=X_{E(h)}=i_{d(E(h))}\zL+E(h)E$). In \cite{VO},
V. Ovsienko proved the noncanonical decomposition
\be\op{Vect}(M)\simeq\op{TVect}(M)\oplus\op{CVect(M)}\label{OvsDecomp}\ee
of $\op{Vect}(M)$ into a direct sum of $\op{CVect}(M)$-modules. An
extension of this decomposition, see \cite{FMP}, will be exploited
below.

\section{Infinitesimal projective contact transformations}\label{InfProjContTrans}

Let us first recall that the symplectic algebra $\op{sp}(2n,\C)$
is the Lie subalgebra of $\op{gl}(2n,\C)$ made up by those
matrices $S$ that verify $JS+\widetilde{S}J=0$, where $J$ is the
symplectic unit. This condition exactly means that the symplectic
form defined by $J$ is invariant under the action of $S$. Since
\[\op{sp}(2n,\C)=\left\{
\left(\begin{array}{cc}
A&B\\
C&D\end{array}\right) : A,B,C,D\in
\op{gl}(n,\C),\tilde{B}=B,\tilde{C}=C,D=-\tilde{A} \right\},\] it
is obvious that
\begin{equation}\label{basis}
\epsilon^j\otimes e_{i+n}+\epsilon^i\otimes e_{j+n}\; (i\le j\le
n),\, -\epsilon^{j+n}\otimes e_i-\epsilon^{i+n}\otimes e_j\; (i\le
j\le n),\, -\epsilon^j\otimes e_i+\epsilon^{i+n}\otimes e_{j+n}\;
(i,j\le n),\end{equation} is a basis of $\op{sp}(2n,\C)$. As
usual, we denote by $(e_1,\ldots,e_{2n})$ the canonical basis of
$\C^{2n}$ and by $(\zve^1,\ldots,\zve^{2n})$ its dual basis. \\

Observe now that the Jacobi (or [first] Lagrange) bracket on a
contact manifold $(M,\za)$ can be built out of contact form $\za$,
see Equation (\ref{JacBrack}), or---in view of the aforementioned
$1$-to-$1$ correspondence---out of the homogeneous symplectic
structure of the symplectization. In the following, we briefly
recall the construction via symplectization, see \cite{PM}, of the
Lagrange bracket, contact vector fields, and the Lie algebra
isomorphism $X:\Ci(M)\raa\op{CVect}(M)$. We then use isomorphism
$X$ to depict Lie subalgebras of contact fields, which play a
central role in this work. As part of our construction is purely
local, we confine
ourselves to the Euclidean setting.\\

Take contact manifold $(\R^{2n+1},\za)$, $\za=i^*\zs$, and its
symplectization $(\zp:\R^{2n+2}_+\raa\R^{2n+1},\zw,\zD)$, see
Example 2. Let ${\cal H}^{\zl}={\cal
H}^\zl(\R^{2n+2}_+)=\{H\in\Ci(\R^{2n+2}_+):\zD H=\zl H\}$ be the
space of homogeneous functions of degree $\zl$. Since $\zw$ has
degree $1$, its contravariant counterpart $\zP$ has degree $-1$,
and the corresponding Poisson bracket verifies $\{{\cal
H}^{\zl},{\cal H}^{\zm}\}_{\zP}\subset{\cal H}^{\zl+\zm-1}$. % It suffices to compute the Lie derivative of the bracket with respect to $\zD$.
In particular, ${\cal H}^1$ is a Lie subalgebra.

If ${\cal X}_H=\{H,.\}_{\zP}$ is the Hamiltonian vector field of a
function $H\in{\cal H}^1$, we
have $[\zD,{\cal X}_H]=0.$ % Simple Lie derivative computation.
Hence, ${\cal X}_H$ is projectable, i.e. $\zp_*{\cal X}_H$ is well-defined. % This means that $\zp_{*(\tau x',\tau)}X_{F,(\tau x',\tau)}$ is independent on $\tau$, see my formulae in my decomposition notes.
As ${\cal X}_H$ is a symplectic vector field, $\zp_*{\cal X}_H$ is
a contact field. The correspondence $\zp_*\circ {\cal X}:{\cal
H}^1(\R^{2n+2}_+)\raa \op{CVect}(\R^{2n+1})$ is obviously a Lie
algebra homomorphism. Remark now that a homogeneous function is
known on the entire fiber $\zp^{-1}(x')=\tau(x',1)$ if it is
specified on the point $(x',1)$. Hence, homogeneous functions are
in fact functions on the base. The correspondence is (for
functions of degree 1) of course \be\chi:{\cal
H}^{1}(\R^{2n+2}_+)\ni H\raa h\in\Ci(\R^{2n+1}),\label{Chi}\ee
with $h(x')=H(x',1)$ and
$H(x',\tau)=\tau^{2}H(\tau^{-1}x',1)=\tau^{2}h(\zp(x',\tau))$
(note that $\zD=\frac{1}{2}{\cal E}$). As every contact vector
field is characterized by a unique base function
$h\in\Ci(\R^{2n+1})$, see Equation (\ref{ContHamMap}), hence by a
unique homogeneous function $H\in{\cal H}^1(\R^{2n+2}_+)$,
morphism $\zp_*\circ {\cal X}$ is actually a Lie algebra
isomorphism.

Map $\chi$ is a vector space isomorphism that allows to push the
Poisson bracket $\{.,.\}_{\zP}$ to the base. The resultant bracket
$\{h,g\}=\chi\{\chi^{-1} h,\chi^{-1} g\}_{\zP}$ is the Lagrange
bracket. Eventually, $\chi$ is a Lie algebra isomorphism.
\begin{center}
\begin{picture}(150,100)\put(50,90){${\cal H}^1(\R^{2n+2}_+)$}
\put(60,85){\vector(-1,-2){30}} \put(90,85){\vector(1,-2){30}}
\put(0,10){$\Ci(\R^{2n+1})$}\put(25,60){$\chi$}\put(120,60){$\zp_*\circ{\cal
X}$} \put(55,12.5){\vector(1,0){40}}\put(72.5,2.5){X}
\put(100,10){$\op{CVect}(\R^{2n+1})$}\end{picture}
\end{center}

It is now easily checked that ``contact Hamiltonian isomorphism''
$X$ is, for any $h\in\Ci(\R^{2n+1})$, given by \be X_h=\zp_*{\cal
X}_{\chi^{-1}h}=\sum_k(\p_{p_k}h\p_{q^k}-\p_{q^k}h\p_{p_k})+{\cal
E}_sh\p_t-\p_th{\cal E}_s-2h\p_t,\label{ContHamMapExpl}\ee where
$(p_1,\ldots,p_n,q^1,\ldots,q^n,t)=(p,q,t)=x'$ are canonical
coordinates in $\R^{2n+1}$ and where ${\cal E}_s=\sum_k
(p_k\p_{p_k}+q^k\p_{q^k})$ is the spatial Euler field. When
comparing this upshot with Equations (\ref{JacHamMap}) and
(\ref{JacBrack}), we get the explicit local form of the Lagrange
bracket.\\

We now depict the aforementioned Lie subalgebras of contact vector
fields as algebras of contact Hamiltonian vector fields of Lie
subalgebras of functions. The algebra of Hamiltonian vector fields
of the Lie subalgebra $\op{Pol}(\R^{2n+1})\subset\Ci(\R^{2n+1})$
of polynomial functions is the Lie subalgebra
$\op{CVect}^*(\R^{2n+1})$ of polynomial contact vector fields,
i.e. contact vector fields with polynomial coefficients. The space
of polynomials admits the decomposition
$\op{Pol}(\R^{2n+1})=\oplus_{r\in\N}\oplus_{k=0}^r{\cal P}^{rk}$,
where ${\cal P}^{rk}$ is the space of polynomials
$t^kP_{r-k}(p,q)$ of homogeneous total degree $r$ that have
homogeneous degree $k$ in $t$. The Lie subalgebra $$\op{Pol}^{\le
2}(\R^{2n+1})=\oplus_{r\le 2}\oplus_{k=0}^r{\cal P}^{rk},$$ which
corresponds via $\chi^{-1}$ to the Lie subalgebra
$$\op{Pol}^2(\R^{2n+2}_+)\subset{\cal H}^1(\R^{2n+2}_+),$$ deserves
particular attention (note that if we set ${\frak g}_{-2}={\cal
P}^{00}, {\frak g}_{-1}={\cal P}^{10}, {\frak g}_0={\cal
P}^{11}\oplus{\cal P}^{20}, {\frak g}_1={\cal P}^{21}$, and
${\frak g}_2={\cal P}^{22}$, we obtain a grading of $\op{Pol}^{\le
2}(\R^{2n+1})$ that is compatible with the Lie bracket [when read
on the symplectic level, this new grading means that we assign the
degree $-1$ to coordinate $\tau$, degree $1$ to $t$, and degree
$0$ to any other coordinate]).

Remember now the Lie algebra isomorphism $J:\op{gl}(m,\R)\ni A\raa
A^{i}_jx^j\p_{x^{i}}\in\op{Vect}^0(\R^m)$ between the algebras of
matrices and of linear vector fields (shifted degree). The inverse
of this isomorphism is the Jacobian map
$J^{-1}:\op{Vect}^0(\R^m)\ni X\raa \p_xX\in\op{gl}(m,\R)$. The Lie
subalgebra $\op{Pol}^2(\R^{2n})={\cal P}^{20}$ is mapped by Lie
algebra isomorphism $X$ (remark that on the considered subalgebra
$X={\cal X}$ and $\{.,.\}=\{.,.\}_{\zP}$) onto a Lie subalgebra of
$\op{CVect}(\R^{2n+1})$ and of $\op{Vect}^0(\R^{2n})$, which in
turn corresponds through Lie algebra isomorphism $J^{-1}$ to a Lie
subalgebra of $\op{gl}(2n,\R)$. A simple computation shows that
the natural basis \be p_ip_j\;(i\le j\le n),\,q^{i}q^j\;(i\le j\le
n),\, p_jq^{i}\;(i,j\le n)\label{BasisPoly}\ee of
$\op{Pol}^2(\R^{2n})$ is transformed by morphism $J^{-1}\circ
{\cal X}$ into the above described basis of $\op{sp}(2n,\R)$, see
Equation (\ref{basis}). It is now clear that $J^{-1}\circ{\cal X}$
is a Lie algebra isomorphism between
$(\op{Pol}^2(\R^{2n}),\{.,.\}_{\zP})$ and
$(\op{sp}(2n,\R),[.,.]_{\circ})$, where $[.,.]_{\circ}$ is the
commutator. We denote by $\op{sp}_{2n}$ the Lie subalgebra of
contact vector fields isomorphic to $\op{Pol}^2(\R^{2n})\simeq
\op{sp}(2n,\R)$.

Eventually, we have the following diagram of Lie algebra
isomorphisms:\\

\begin{center}
\begin{picture}(150,100)\put(0,90){$\op{Pol}^2(\R^{2n+2}_+)$}
\put(30,85){\vector(0,-2){60}} \put(35,85){\vector(2,-1){120}}
\put(0,10){$\op{Pol}^{\le
2}(\R^{2n+1})$}\put(15,60){$\chi$}\put(90,60){$\zp_*\circ{\cal
X}$}
\put(65,12.5){\vector(1,0){70}}\put(160,85){\vector(0,-2){60}}\put(65,92,5){\vector(1,0){70}}
\put(80,97.5){$J^{-1}\circ{\cal X}$}\put(95,2.5){X}
\put(150,10){$\op{sp}_{2n+2}$}\put(140,90){$\op{sp}(2n+2,\R)$}\end{picture}
\end{center}

\noindent It is obvious that the right bottom algebra is a Lie
subalgebra of contact vector fields that is isomorphic with
$\op{sp}(2n+2,\R)$: hence the notation. The right vertical arrow
refers to the embedding of $\op{sp}(2n+2,\R)$ into
$\op{CVect}(\R^{2n+1})$ that can be realized just as the
projective embedding of $\op{sl}(m+1,\R)$ into $\op{Vect}(\R^m)$.
More precisely, the linear symplectic group $\op{SP}(2n+2,\R)$
naturally acts on $\R^{2n+2}$ by linear symplectomorphisms. The
projection $\rho(S)(x') := \pi(S.(x',1)),$ $S\in\op{SP}(2n+2,\R)$,
$x'\in\R^{2n+1}$, of this action ``$.$'' induces a ``local''
action on $\R^{2n+1}$. The tangent action to projection $\zr$ is a
Lie algebra homomorphism that maps the symplectic algebra
$\op{sp}(2n+2,\R)$ into contact vector fields
$\op{CVect}(\R^{2n+1})$. We refer to the Lie subalgebra generated
by the fundamental vector fields associated with this
infinitesimal action as the algebra of \emph{infinitesimal
projective contact transformations}. This algebra $\op{sp}_{2n+2}$
is a maximal proper Lie subalgebra of $\op{CVect}^*(\R^{2n+1})$
(just as the projective embedding $\op{sl_{m+1}}$ of
$\op{sl}(m+1,\R)$ is a maximal proper Lie subalgebra of
$\op{Vect}^*(\R^m)$). Over a Darboux chart, any
$(2n+1)$-dimensional contact manifold can be identified with
$(\R^{2n+1},i^*\zs)$. It is therefore natural to consider
$\op{sp}_{2n+2}$ as a subalgebra of vector fields over the
chart.\\

Eventually, a basis of $\op{sp}_{2n+2}$ can be deduced via
isomorphism $X$ from the canonical basis of $$\op{Pol}^{\le
2}(\R^{2n+1})=\oplus_{r\le 2}\oplus_{k=0}^r{\cal P}^{rk}.$$ Using
Equation (\ref{ContHamMapExpl}), we immediately verify that the
contact Hamiltonian vector fields of $1\in{\cal P}^{00}$ and
$p_i,q^{i}\in {\cal P}^{10}$ are \be X_1=-2\p_{t}=E,
X_{p_i}=\p_{q^{i}} - p_i\p_{t}, X_{q^{i}}=-\p_{p_i} -
q^{i}\p_{t}.\label{BasisSp1}\ee These fields generate a Lie
algebra $\tilde{{\frak h}}_{n,1}$ that is isomorphic to the
Heisenberg algebra ${\frak h}_n$. Let us recall that the
Heisenberg algebra ${\frak h}_n$ is a nilpotent Lie algebra with
basis vectors $(a_1,\ldots,a_n,b_1,\ldots,b_n,c)$ that verify the
commutation relations
$$[a_i,b_j]=\zd_{ij}c, [a_i,a_j]=[b_i,b_j]=[a_i,c]=[b_i,c]=0.$$
Similarly the Hamiltonian vector field of $t\in{\cal P}^{11}$ is
the modified Euler field \be X_t=-{\cal
E}_s-2t\p_t,\label{BasisSp2}\ee and the Hamiltonian vector fields
of $p_ip_j, q^{i}q^j,p_jq^{i}\in{\cal P}^{20}$ are \be
X_{p_ip_j}=p_j\p_{q^{i}}+p_i\p_{q^j}\;(i\le j\le n),\,
X_{q^{i}q^j}= -q^j\p_{p_i}-q^{i}\p_{p_j}\;(i\le j\le n),\,
X_{p_jq^{i}}=q^{i}\p_{q^j}-p_j\p_{p_i}\;(i,j\le
n).\label{BasisSp3}\ee These fields form the basis of
$\op{sp}_{2n}$ that corresponds via $J$ to the basis of
$\op{sp}(2n,\R)$ specified in Equation (\ref{basis}). Finally,
$p_it, q^{i}t\in{\cal P}^{21}$, and $t^2\in{\cal P}^{22}$ induce
the fields \be X_{p_it}=t(\p_{q^{i}}-p_i\p_t)-p_i{\cal
E}_s=t\p_{q^{i}}-p_i{\cal E},\,
X_{q^{i}t}=t(-\p_{p_i}-q^{i}\p_t)-q^{i}{\cal
E}_s=-t\p_{p_{i}}-q^{i}{\cal E},\, X_{t^2}=-2t{\cal
E}.\label{BasisSp4}\ee Again these fields generate a Lie algebra
$\tilde{{\frak h}}_{n,2}$ that is a model of the Heisenberg
algebra ${\frak h}_n$.

Observe also that the contact Hamiltonian vector field of a member
of ${\cal P}^{rk}$ ($r\ge 3$) is a polynomial contact field of
degree $r$. Finally, the algebra $\tilde{{\frak h}}_{n,1}\oplus\R
X_t\oplus\op{sp}_{2n}$ is the algebra
$\op{AVect}(\R^{2n+1})\cap\op{CVect}(\R^{2n+1})$ of affine contact
vector fields.

\section{Differential operators, symbols, actions, tensor
densities}\label{DiffOpSymbEct}

Let $\zp:E\raa M$ and $\zt:F\raa M$ be two (finite rank) vector
bundles over a (smooth $m$-dimensional) manifold $M$.\\

We denote by ${\cal D}_k(E,F)$, $k\in\N$, the space of {\it $k$th
order linear differential operators} between the spaces
$\zG^{\infty}(E)$ and $\zG^{\infty}(F)$ of smooth global sections
of $E$ and $F$ (in the following we simply write $\zG(E)$ or
$\zG(F)$), i.e. the space of the linear maps
$D\in\op{Hom}_{\R}(\zG(E),\zG(F))$ that factor through the $k$th
jet bundle $J^kE$ (i.e. for which there is a bundle map
$\tilde{D}:J^kE\raa F$, such that $D=\tilde{D}\circ i^0$, where
$i^0:E\raa J^kE$ is the canonical injection and where the {\small
RHS} is viewed as a map between sections). %An equivalent purely
% algebraic definition ``\`a la Vinogradov'' has been used in
% \cite{GP}: $k$th order linear differential operators are exactly
% the linear maps $D\in\op{Hom}_{\R}(\zG(E),\zG(F))$ that verify
% $[\ldots [D,f_0]\ldots f_{k}]=0$, for any
% $f_0,\ldots,f_k\in\Ci(M)$ (here $[D,f](s):=D(fs)-fD(s)$,
% $s\in\zG(E)$).
It is obvious that $0$ -- order differential
operators are just the sections $\zG(\op{Hom(E,F)})\simeq
\zG(E^*\otimes F)$ and that the space ${\cal D}(E,F)=\cup_k{\cal
D}_k(E,F)$ of all linear differential operators between $E$ and
$F$ (or better between $\zG(E)$ and $\zG(F)$) is filtered by the
order of differentiation.

The {\it $k$th order principal symbol} $\zs_k(D)$ of an operator
$D\in{\cal D}_k(E,F)$ is the map $\tilde{D}\circ i^k$, where $i^k$
denotes the canonical injection $i^k:{\cal S}^kT^*M\otimes E\raa
J^kE$. Actually, this compound map is a bundle morphism $\zs_k(D):
{\cal S}^kT^*M\otimes E\raa F$, or, equivalently, a section
$\zs_k(D)\in\zG({\cal S}^kTM\otimes E^*\otimes F)$. % The following
% alternative viewpoint favors the algebraic approach. Any $k$th
% order operator $D\in{\cal D}_k(E,F)$ defines a map
% $\tilde{\zs}_k(D):\Ci(M)^{\times k}\ni(f_1,\ldots,f_k)\raa
% [\ldots[D,f_1]\ldots f_k]\in{\cal D}_0(E,F)$. As
% $\tilde{\zs}_k(D)$ only depends on $(df_1,\ldots,df_k)$ and is
% symmetric in its arguments, it implements a symmetric
% $\Ci(M)$-multilinear map $\zs_k(D)$ between $\zG(T^*M)^{\times k}$
% and $\zG(E^*\otimes F)$. This map is of course a section
% $\zs_k(D)\in\zG({\cal S}^kTM\otimes E^*\otimes F)$.
In the following, we call symbol space (associated with ${\cal
D}(E,F)$), and denote by ${\cal S}(E,F)$, the graded space ${\cal
S}(E,F)=\oplus_k{\cal S}^k(E,F)$, where ${\cal
S}^k(E,F):=\zG({\cal S}^kTM\otimes E^*\otimes F)$. Since
$\zs_k:{\cal D}_k(E,F)\raa {\cal S}^k(E,F)$ is a linear
surjection, it induces a vector space isomorphism between the
graded space associated with the filtered space ${\cal D}(E,F)$
and the graded space ${\cal S}(E,F)$.\\ % Of course, symbols can be
% interpreted as polynomials along the fibers of the cotangent
% bundle $T^*M$, with coefficients in the linear maps between the
% fibers of $E$ and $F$, which depend smoothly on $x$ in
% the base $M$.\\
% Note that ${\cal S}^kT^*M$ are the $k$th order derivatives and that $E$ are the components.
% Jets: $\p^{\za}_x s^j$, $\mid\za\mid \le r$, $x\in M$, $s$ section defined around $x$, $s^j$ component

Roughly spoken, an {\it equivariant or natural quantization} is a
vector space isomorphism $Q:{\cal S}(E,F)\raa {\cal D}(E,F)$ that
verifies some normalization condition and intertwines the actions
on ${\cal S}(E,F)$ and ${\cal D}(E,F)$ of some symmetry group $G$
of base manifold $M$. However, in order to define such actions,
the action $\zf^M$ of $G$ on $M$ should lift to $E$ (and $F$) as
an action $\zf^{E}$ (resp. $\zf^{F}$) of $G$ by vector bundle maps
$\zf^{E}_g:E\raa E$, $g\in G$, over the corresponding maps
$\zf^M_g:M\raa M$. Actually, the action $\zf^{\zG(E)}$ of $G$ on
$\zG(E)$  can then be defined by
$$(\zf^{\zG(E)}_gs)_x:=\zf^{E}_{g^{-1}}s_{\zf^M_g(x)},$$ $g\in G$,
$s\in\zG(E)$, $x\in M$, and the action $\zf^{\cal D}$ of $G$ on
${\cal D}(E,F)$ is $$\zf^{\cal D}_gD:=\zf^{\zG(F)}_{g}\circ
D\circ\zf^{\zG(E)}_{g^{-1}},$$ for any $g\in G, D\in{\cal
D}(E,F)$. Eventually, there is also a canonical action $\zf^{\cal
S}$ on symbols. Indeed, for any $g\in G$, $P\in{\cal
S}^k(E,F)=\zG({\cal S}^kTM\otimes E^*\otimes F)$, $x\in M$, $e\in
E_x$, it suffices to set
$$(\zf^{\cal S}_gP)_x(e):=({\cal S}^kT\zf^M_{g^{-1}}\otimes\zf^{F}_{g^{-1}})P_{\zf^M_g(x)}(\zf^{E}_ge).$$

The appropriate setting for such investigations is the framework
of natural functors (for all questions related with natural
functors and natural operations, we refer the reader to
\cite{KMS}, for a functorial approach to natural quantization, see
\cite{MB}). % (we confine ourself here to recalling that a natural
% fiber [resp. vector] bundle functor is a covariant functor from
% the category of $m$-dimensional smooth manifolds and local
% diffeomorphisms into the category of fiber [resp. vector] bundles
% and bundle maps, which verifies [natural] base and locality
% conditions).
Indeed, let $\mathbb{F}$ and $\mathbb{F'}$ be two natural vector
bundle functors and consider differential operators and symbols
between the vector bundles $E=\mathbb{F}M$ and $F=\mathbb{F'}M$
over an $m$-dimensional smooth manifold $M$. If now $\zf^M_g:M\raa
M$ is a local diffeomorphism, then
$\zf^{E}_g=\mathbb{F}\zf^M_g:\mathbb{F}M\raa\mathbb{F}M$ (resp.
$\zf^F_g=\mathbb{F'}\zf^M_g$) is a vector bundle map over
$\zf^M_g$. Hence, actions on the base lift canonically and actions
of the group $\op{Diff}(M)$ of local diffeomorphisms of $M$ (and
of the algebra $\op{Vect}(M)$ of vector fields of $M$) can be
defined (as detailed above) on sections $\zG(\mathbb{F}M)$ and
$\zG(\mathbb{F}'M)$, as well as on differential operators ${\cal
D}(\mathbb{F}M,\mathbb{F}'M)$ and symbols ${\cal
S}(\mathbb{F}M,\mathbb{F}'M)$ between these spaces of sections.

Remember now that there is a $1$-to-$1$ correspondence between
representations of the jet group $G^r_m$ on vector spaces $V$ and
natural vector bundle functors $\mathbb{F}$ of order $r$ on the
category of $m$-dimensional smooth manifolds $M$, see
\cite[Proposition 14.8]{KMS}. The objects of such a functor are
the vector bundles $\mathbb{F}M=P^rM\times_{G^r_m}V$ associated
with the $r$th order frame bundles $P^rM$. So the canonical
representation of $G^1_m=\op{GL}(m,\R)$ on the (rank $1$) vector
spaces $\zD^{\zl}\R^m$ ($\zl\in\R$) of $\zl$-densities on $\R^m$
induces a $1$-parameter family of natural $1$st order vector
bundle functors $\mathbb{F}_{\zl}$. Hence, we get $\op{Diff}(M)$-
and $\op{Vect}(M)$-actions on sections of the (trivial) line
bundles
$$\mathbb{F}_{\zl}M=P^1M\times_{\op{GL}(m,\R)}\zD^{\zl}\R^m=\zD^{\zl}TM$$
of $\zl$-densities of $M$, i.e. on tensor densities ${\cal
F}_{\zl}(M):=\zG(\mathbb{F}_{\zl}M)$ of order $\zl$ of $M$. As
aforementioned these actions generate actions on differential
operators ${\cal D}_{\zl\zm}(M):={\cal
D}(\mathbb{F}_{\zl}M,\mathbb{F}_{\zm}M)$ between tensor densities
of weights $\zl$ and $\zm$, and on the corresponding symbols
${\cal S}_{\zd}(M):={\cal S}(\mathbb{F}_{\zl}M,\mathbb{F}_{\zm}M)$
$=\zG({\cal S}TM\otimes
\mathbb{F}_{\zl}^*M\otimes\mathbb{F}_{\zm}M)=\zG({\cal S}TM\otimes
\mathbb{F}_{\zd}M),$ where $\zd=\zm-\zl.$

Primarily the local forms of these actions are well-known. Below,
we focus on the algebra actions rather than on the group actions.
Let us recall that triviality of the line bundles
$\mathbb{F}_{\zl}M$ has been proven via construction of a nowhere
vanishing section $\zr_0$ of $\mathbb{F}_1M=\zD^1TM$ that has at
each point only strictly positive values. If the considered
manifold $M$ is orientable, we can set $\zr_0=\mid\!\zW\!\mid$,
where $\zW$ is a volume of $M$. Let us choose such a
trivialization $\zr_0$. The correspondences $\zt_0^{\zl}:\Ci(M)\ni
f\raa f\zr_0^{\zl}\in{\cal F}_{\zl}(M)$ are then vector space
isomorphisms and the actions $L^{\zl}$ of vector fields on the
spaces ${\cal F}_{\zl}(M)$ are, for any $X\in\op{Vect}(M)$ and any
$f\in\Ci(M)$, given by
$$L^{\zl}_X(f\zr_0^{\zl})=(X(f)+\zl f\op{div}_{\zr_0}X)\zr^{\zl}_0.$$
Furthermore, for any $X\in\op{Vect}(M)$, $D\in{\cal
D}_{\zl\zm}(M)$, and $P\in{\cal S}_{\zd}(M)$, we have
$${\cal L}_XD=L^{\zm}_X\circ D-D\circ L^{\zl}_X,$$
and \be L_XP=X^{\sharp}(P)+\zd P\op{div}_{\zr_0}
X,\label{ActionSymb}\ee where $X^{\sharp}$ denotes the cotangent
lift of $X$ and where we have omitted in the {\small LHS} the
dependance of the actions on $\zl$, $\zm$, and on $\zd$,
respectively.
% If we choose a trivialization $\zr_0$ of the $1$-density bundle, we find that
% $\zG(E\otimes\zD^{\zl}TM) \simeq \zG(E)\zr_0^{\zl}.$

\section{Invariant tensor fields}\label{InvTensFields}

Consider a $(2n+1)$-dimensional smooth Hausdorff second countable
(coorientable) contact manifold $(M,\za)$. Let us recall that this
work is originated from the classification problem of the spaces
$({\cal D}_{\zl\zm}(M),{\cal L})$ as modules over the Lie algebra
$\op{CVect}(M)$ of contact vector fields. A first approximation is
the computation of the intertwining operators $T$ between the
corresponding $\op{CVect}(M)$-modules $({\cal S}_{\zd}(M),L)$.
Note that locally these symbol spaces are also modules over the
Lie subalgebra $\op{sp}_{2n+2}\subset\op{CVect}(M)$. If such a
module morphism
$$T:{\cal S}_{\zd}^{\E}(M)=\zG({\cal
S}^{\E}TM\otimes\mathbb{F}_{\zd}M)\raa {\cal
S}_{\ze}^m(M)=\zG({\cal S}^mTM\otimes\mathbb{F}_{\ze}M)$$ is a
$k$th order differential operator, its principal symbol
$\zs_k(T)\in\zG({\cal S}^kTM\otimes{\cal S}^mTM\otimes{\cal
S}^{\E}T^*M\otimes\mathbb{F}_{\zn}M),$ $\zn=\ze-\zd$, is (roughly
spoken) again invariant, see below. Hence, the quest for tensor
fields in the preceding symbol space \be {\cal
S}^{km}_{\E;\zn}(M):=\zG({\cal S}^kTM\otimes{\cal
S}^mTM\otimes{\cal
S}^{\E}T^*M\otimes\mathbb{F}_{\zn}M),\label{SymbolSpace}\ee
$k,m,\E\in\N$, $\zn\in\R$, which are $\op{CVect}(M)$- and,
locally, $\op{sp}_{2n+2}$-invariant for the canonical action.

\subsection{$\mathbf{\op{{\bf CVect}}(M)}$- and $\mathbf{\op{\bf sp}_{2n+2}}$-invariants}

In the following, we need a result on Taylor expansions. If $f\in
C^{\infty}(\R^m)$ and $x_{0}\in\R^m$, we denote by
$t^k_{x_{0}}(f)$ the $k$th order Taylor expansion of $f$ at
$x_{0}$. We use the same notation for Taylor expansions of vector
fields.
\begin{prop}\label{taylor}For every $f\in C^{\infty}(\R^{2n+1})$ and $x_{0}\in \R^{2n+1}$,
we have
        \[t_{x_{0}}^1(X_{f}) =t_{x_{0}}^1(X_{t^2_{x_{0}}(f)}).\]
        \end{prop}
        {\it Proof}.
            First note that, in any coordinate system
            $(x^1,\ldots,x^{2n+1})$, if ${\cal E}= x^i\p_{x^{i}}$ is the
            Euler field, one has, for all $k\geq 1$,
            \[\left\{\begin{array}{lll}
        \p_{x^{i}}\circ t^k_{x_{0}} & = &t^{k-1}_{x_{0}}\circ\p_{x^{i}},\\
        ({\cal E}-k\,\mathrm{id}) \circ t^{k}_{x_{0}} & = &
        t^{k-1}_{x_{0}}\circ({\cal E}-k\,\mathrm{id}).
        \end{array}\right.\]
        This allows to show that
        \[t^1_{x_{0}}(X_{f}) = X_{t^2_{x_{0}}(f)} +
        [t^1_{x_{0}}(\p_{x^{2n+1}}f)-(\p_{x^{2n+1}}f)(x_{0})][{\cal E}-{\cal E}_{x_{0}}].\]
        The result follows, since the last term of the {\small RHS} and its partial derivatives
        vanish at $x_{0}$.\quad\rule{1.5mm}{2.5mm}\\

It is clear that the spaces
$\zG(\otimes^p_{q}TM\otimes\mathbb{F}_{\zn}M)$ are again
representations of the Lie algebra $\mathrm{CVect}(M)$.

\begin{theo}\label{sp}
    Let $M$ be a (coorientable) contact manifold of dimension $2n+1$. A tensor field
    $u\in \Gamma(\otimes^p_{q}TM\otimes\mathbb{F}_{\zn}M)$ is $\mathrm{CVect}(M)$-invariant if and only if, over any
    Darboux chart, $u$ is $\mathrm{sp}_{2n+2}$-invariant.
    \end{theo}
{\it Proof}.
      The Lie derivative
      \[L : \mathrm{Vect}(M)\times
      \Gamma(\otimes^p_{q}TM\otimes\mathbb{F}_{\zn}M)\raa
      \Gamma(\otimes^p_{q}TM\otimes\mathbb{F}_{\zn}M) :
      (X,u)\raa L_{X}u\]
      is a differential operator that has order $1$ in the first argument.
      In other words, the value $(L_{X}u)_{x_{0}}$, $x_0\in M$, only depends on the
      first jet $j^1_{x_0}(X)$ of $X$ at $x_{0}$. Hence, the result follows from
      Proposition \ref{taylor}.\quad\rule{1.5mm}{2.5mm}

\subsection{Particular invariants}

We continue to work on a ($2n+1$)-dimensional coorientable contact
manifold $(M,\za)$ endowed with a fixed contact form, and describe
basic contact invariant or locally affine contact invariant tensor
fields in ${\cal S}(M):=\oplus_{km\E\zn}{\cal
S}^{km}_{\E;\zn}(M)$, see Equation (\ref{SymbolSpace}). Observe
that locally we can view the elements of ${\cal
S}^{km}_{\E;\zn}(M)$ as polynomials of homogeneous degrees $k,$
$m$, and $\E$ in fiber variables $\xi,$ $\zh,$ and $Y$, with
coefficients in ${\cal F}_{\zn}(M)$.

\begin{enumerate}
\item The {\it identity endomorphism} $u\in\zG(TM\otimes T^*M)$ of
$TM$ can be viewed as a $\op{Vect}(M)$-invariant element
$u_1\in{\cal S}^{10}_{1;0}(M)$ and as a $\op{Vect}(M)$-invariant
element $u_2\in{\cal S}^{01}_{1;0}(M)$. Locally, these two
invariant tensor fields read $u_{1}:(\xi,\eta,Y)\raa\langle
Y,\xi\rangle$ and $u_{2}:(\xi,\eta,Y)\raa\langle Y,\eta\rangle,$
where $\la.,.\ra$ denotes the contraction.

\item {\it Contact form} $\alpha$ induces a
$\op{CVect}(M)$-invariant tensor field $$u_3=\za\,\otimes
\vert\zW\vert^{-I}\in{\cal S}^{00}_{1;-\op{I}}(M),$$
$\zW=\za\w(d\za)^n$, $\op{I}=\frac{1}{n+1}.$ Invariance of $u_3$
with respect to the action of contact fields is a direct
consequence of Equation (\ref{DivContact}). The local form of
$u_3$ is $u_{3}:(\xi,\eta,Y)\raa
\alpha(Y)\vert\Omega\vert^{-\op{I}}.$

\item The next invariant tensor field is implemented by the {\it
Lagrange bracket}. Let us first mention that our construction of
the Lagrange bracket on $(\R^{2n+1},i^*\zs)$, as pullback of the
Poisson bracket of the symplectization of this contact structure,
can be generalized to an arbitrary contact manifold $M$, see
\cite[Section 10.2, Corollary 2]{OR}: the Poisson bracket on the
symplectization of $M$ defines a bracket $\{.,.\}$, called
Lagrange bracket, on the space ${\cal
F}_{-\op{I}}(M)=\zG(\mathbb{F}_{-\op{I}}M)$ of $(-\op{I})$-tensor
densities of $M$. This bracket is a first order bidifferential
operator between ${\cal F}_{-\op{I}}(M)\times {\cal
F}_{-\op{I}}(M)$ and ${\cal F}_{-\op{I}}(M)$. Its principal symbol
is defined just along the same lines than the principal symbol of
a differential operator, see Section \ref{DiffOpSymbEct}. Hence,
this symbol $\zs_{11}(\{.,.\})$ is a tensor field
$$L_1:=\zs_{11}(\{.,.\})\in\zG(TM\otimes
TM\otimes\mathbb{F}_{\op{I}}M)={\cal S}^{11}_{0;\op{I}}(M).$$ It
is a basic fact in equivariant quantization that the principal
symbol of a (multi)differential operator between tensor densities
intertwines the actions $L$ and ${\cal L}$ of vector fields on
symbols and operators (a short computation in local coordinates
also allows to assure oneself of this fact).
$\op{CVect}(M)$-invariance of $L_1$ then follows from contact
invariance of $\{.,.\}$ that in turn is nothing but a
reformulation of the Jacobi identity, see Equation
(\ref{LieDensities}). The local polynomial form of $L_1$ follows
from Equation (\ref{ContHamMapExpl}): $L_{1}:(\xi,\eta,Y)\raa
(\sum_k(\xi_{p_k}\zh_{q^k}-\xi_{q^k}\zh_{p_k})+\zh_t\la{\cal
E}_s,\xi_s\ra-\xi_t\la{\cal E}_s,\zh_s\ra)\vert\zW\vert^{\op{I}}$,
where $\zb_s=\zb-\zb_tdt$, for any one form $\zb$.

\item Eventually, the {\it Reeb vector field} $E$ also induces an
invariant tensor field
$u=E\otimes\vert\zW\vert^{\op{I}}\in\zG(TM\otimes\mathbb{F}_{\op{I}}M)$,
which can be viewed as an element $u_4\in{\cal
S}^{10}_{0;\op{I}}(M)$ and as an element $u_5\in{\cal
S}^{01}_{0;\op{I}}(M)$. Since $[X_f,E]=-i_{d(E(f))}\zL-E(f)E$
(confer end Section \ref{ContGeo}), it is easily seen that fields
$u_4$ and $u_5$ are not contact invariant, but only (locally)
affine contact invariant (confer end Section
\ref{InfProjContTrans}). They (locally) read $u_{4}:
(\xi,\eta,Y)\raa -2\xi_{t}\vert\Omega\vert^{\op{I}}$ and
$u_{5}:(\xi,\eta,Y)\raa -2\zh_t\vert\zW\vert^{\op{I}}.$
\end{enumerate}

\subsection{Classifications}

In this subsection we classify affine contact invariant and
contact invariant tensor fields.

\begin{theo}\label{invaff}
    The polynomials $u_{i},\,i\in\{1,\ldots,5\}$, and $L_{1}$ generate the
    algebra of $\op{AVect}(\R^{2n+1})\cap\op{CVect}(\R^{2n+1})$-invariant polynomials
    in ${\cal S}(\R^{2n+1})$.
    \end{theo}

{\it Proof}. In the following we refer to the algebra of invariant
polynomials generated by $u_{1},\ldots,u_{5},$ and $L_{1}$ as the
space of {\it classical invariant polynomials}. In order to show
that there are no other invariants, we prove that the dimensions
of the subspace ${\cal S}_1$ of classical invariant polynomials in
${\cal S}^{km}_{\E;\zn}(\R^{2n+1})$ and of the subspace ${\cal
S}_2$ of all invariant polynomials inside ${\cal
S}^{km}_{\E;\zn}(\R^{2n+1})$ coincide, for any fixed
$(k,m,\E,\zn)$.\\

       Since the polynomials
       \[ u_{1}^a\,u_{2}^b\,u_{3}^c\,u_{4}^d\,u_{5}^e\,L_{1}^f\in{\cal S}^{a+d+f,b+e+f}_{a+b+c;\op{I}(d+e+f-c)}(\R^{2n+1})\]
       are independent and belong to ${\cal S}^{km}_{\E;\zn}(\R^{2n+1})$
       if and only if $(a,b,c,d,e,f)\in\N^6$ is a solution of
       \[(S_1) :\left\{\begin{array}{lll}
       a+d+f&=&k\\
       b+e+f &=&m\\
       a+b+c & = &\E\\
       d+e+f-c&=&(n+1)\zn
       \end{array}\right.,\]
       the dimension of ${\cal S}_1$
       is exactly the number (which is clearly finite) of solutions in $\N^6$ of system
       $(S_1)$. ($\star$)\\

       Let now $Q\in{\cal S}^{km}_{\E;\zn}(\R^{2n+1})$ be an arbitrary invariant
       polynomial and set
       \[Q(\xi,\eta,Y) = \sum_{i=0}^k\sum_{j=0}^m\sum_{r=0}^{\E}
       \xi_{t}^i\eta_{t}^j(Y^t)^rQ_{i,j,r}(\xi_{s},\eta_{s},Y_{s}),\]
       where polynomial $Q_{i,j,r}$ is homogeneous of degree
       $(k-i,m-j,\E-r)$.
       The degree defined by
       \[\op{Deg}(\xi_{t}^i\eta_{t}^j(Y^t)^rQ_{i,j,r}(\xi_{s},\eta_{s},Y_{s}))
       :=i+j-r\] will be basic in our investigation.
       Let us recall that the obvious extension of action (\ref{ActionSymb}) to ${\cal S}^{km}_{\E;\zn}(\R^{2n+1})$
       reads, for all vector fields $X\in\op{Vect}(\R^{2n+1})$, \be L_XQ=X(Q)-\p_jX^{i}\eta_i\p_{\eta_j}Q-\p_jX^{i}\xi_i\p_{\xi_j}Q
       +\p_jX^{i}Y^j\p_{Y^{i}}Q+\zn\p_iX^{i}Q\label{invcond},\ee where $\p_k$ denotes the derivative with respect to the
       $k$th coordinate of $\R^{2n+1}$. It is easily checked that the invariance conditions
       with respect to the contact Hamiltonian vector fields $X_{1},\,X_{p_{\frak a}},\,X_{q^{\frak a}}$, and $X_{t}$ (${\frak a}\in\{1,\ldots, n\}$),
       see Section \ref{InfProjContTrans}, read
       \[\left\{\begin{array}{lll}
   \p_{t}Q_{i,j,r} & = & 0\\
   \p_{p_{\frak a}} Q_{i,j,r} - \p_{\xi_{q^{\frak a}}} Q_{i-1,j,r} - \p_{\eta_{q^{\frak a}}}
   Q_{i,j-1,r} +(r+1)Y^{q^{\frak a}} Q_{i,j,r+1} & = &0\\
   \p_{q^{\frak a}} Q_{i,j,r} + \p_{\xi_{p_{\frak a}}} Q_{i-1,j,r} + \p_{\eta_{p_{\frak a}}}
   Q_{i,j-1,r} - (r+1)Y^{p_{\frak a}} Q_{i,j,r+1} & = &0\\
   {\cal E}_{s}(Q_{i,j,r}) - [(k+m-\E) +(i+j-r)-2(n+1)\zn]\,Q_{i,j,r} &=&0
   \end{array}\right..\]
   Let $d_0$ be the lowest degree $\op{Deg}$ in $Q$. If $Q_{i,j,r}$ is part of a term of degree $d_0$, the first
   three equations of the above system imply that $Q_{i,j,r}$ has constant
   coefficients. The fourth equation entails that
   $d_{0}=2(n+1)\zn-(k+m-\E).$
   An easy induction shows that all polynomials
   $Q_{i,j,r}$ have polynomial coefficients and that they are
   completely determined by the lowest degree terms. So, the dimension
   of the space ${\cal S}_2$ of invariant polynomials in ${\cal S}^{km}_{\E;\zn}(\R^{2n+1})$ is at most the
   dimension of the space of lowest degree terms. ($\star\star$)

   We now take a closer look at these lowest degree terms
   \begin{equation}\label{low}\sum_{i+j-r =d_{0}}
 \xi_{t}^i\eta_{t}^j(Y^t)^rQ_{i,j,r}(\xi_{s},\eta_{s},Y_{s}),\end{equation}
   where the polynomials $Q_{i,j,r}$ have constant coefficients, and use the invariance conditions with respect to the algebra
   $\mathrm{sp}_{2n}$. Observe first that the Lie derivatives in the direction of
   the fields of this algebra preserve the degree $\op{Deg}$. Indeed, for any field $X$ of the basis of
   $\op{sp}_{2n}$, see Section \ref{InfProjContTrans}, the derivatives $\p_jX^{i}$ vanish for $X^t$ and for $\p_t.$
   Hence, every polynomial $Q_{i,j,r}$ in
   (\ref{low}) must be $\mathrm{sp}_{2n}$-invariant. As these polynomials have constant coefficients and the
   considered vector fields have vanishing divergence, this means that any $Q_{i,j,r}$ in (\ref{low})
   is invariant for the canonical $\op{sp}(2n,\R)$-action. When applying a classical result of Weyl, \cite{HW}, we conclude
   that each polynomial $Q_{i,j,r}$ in (\ref{low}) is a polynomial in the variables $\langle
   Y_{s},\xi_{s}\rangle$, $\langle Y_{s},\eta_{s}\rangle$, and
   $\Pi(\xi_{s},\eta_{s})$, where $\Pi=\sum_{\frak a}\p_{p_{\frak a}}\w\p_{q^{\frak
   a}}$. Eventually, the lowest degree terms of $Q$ read
   \[\sum_{i,j,r,\za,\zb,\zg\in\N}c_{ijr\za\zb\zg}\,\;\xi_{t}^i\eta_{t}^j
   (Y^t)^r \langle Y_{s},\xi_{s}\rangle^{\za} \langle
   Y_{s},\eta_{s}\rangle^{\zb} \Pi(\xi_{s},\eta_{s})^\zg,\]
   where $c_{ijr\za\zb\zg}\in\R$ and where
   $(i,j,r,\za,\zb,\zg)\in\N^6$ is a solution of the system
   \[(S_2):\left\{\begin{array}{lll}
\za+\zg &=&k-i\\
\zb+\zg &=&m-j\\
\za+\zb &=&\E-r\\
i + j -r &= &  2(n+1)\zn -(k+m-\E)
\end{array}\right..\]
   This system implies in particular that $(n+1)\zn=k+m-\E-\zg$ is an
   integer. It is easily checked that system ($S_2$) is equivalent to
   system ($S_1$). When taking into account upshots ($\star$) and ($\star\star$),
   we finally see that the dimension of the space ${\cal S}_2$ of all invariant polynomials in ${\cal S}^{km}_{\E;\zn}(\R^{2n+1})$
   is at most the dimension of the space ${\cal S}_1$ of classical
   invariant polynomials in ${\cal S}^{km}_{\E;\zn}(\R^{2n+1})$.
   \rule{1.5mm}{2.5mm}\\

   As a corollary, we get the following
         \begin{theo}\label{invcont}
    For every (coorientable) contact manifold $M$, the fields $u_{1}, u_{2}, u_{3}$, and $L_{1}$ generate the
    algebra of $\op{CVect}(M)$-invariant fields in ${\cal
    S}(M)$.
    \end{theo}

      {\it Proof}. Any $\op{CVect}(M)$-invariant field in ${\cal S}^{km}_{\E;\zn}(M)$ is over every Darboux
      chart an $\op{sp}_{2n+2}$-invariant polynomial, Theorem \ref{sp}. Hence, due to Theorem
      \ref{invaff}, it reads as a linear combination of
      polynomials $u_1^{a}u_2^bu_3^cu_4^du_5^{e}L_1^f$. Invariance with respect to the second Heisenberg algebra ${\tilde{\frak h}}_{n,2}$, see Section
      \ref{InfProjContTrans}, allows to satisfy oneself that $d=e=0.$ Computations are straightforward (but tedious) and
      will not be given here. \rule{1.5mm}{2.5mm}

\section{Invariant operators between symbol modules implemented by the same density
weight}\label{InvOpSameDens}

We now use Theorem \ref{invcont} concerning contact-invariant
tensor fields to classify specific ``classical'' module morphisms,
i.e. intertwining operators between $\op{sp}_{2n+2}$-modules of
symbols induced by the same density
weight, see below.\\

We first recall the definition of two invariant operators that
were basic in \cite{FMP}.

Let $(M,\za)$ be a Pfaffian manifold. In the following, we denote
the $\op{CVect}(M)$-invariant tensor field
$u_3=\za\otimes\vert\zW\vert^{-\op{I}}\in{\cal
S}^{00}_{1;-\op{I}}(M)=\zG(T^*M\otimes\mathbb{F}_{-\op{I}}M)$
simply by $\za$ (if more precise notation is not required in order
to guard against confusion). Contact form $\alpha$ can then be
viewed as a contraction operator $$i_{\za}:{\mathcal
S}^k_\delta(M)\raa {\mathcal S}^{k-1}_{\delta-\op{I}}(M),$$ where
${\cal S}^k_{\zd}(M)=\zG({\cal S}^kTM\otimes\mathbb{F}_{\zd}M)$.
Due to invariance of $\za$, the vertical cotangent lift $i_{\za}$
of $\za$ is clearly a $\op{CVect}(M)$-intertwining operator.

We also extend the contact Hamiltonian operator, see Equation
(\ref{ContHamMapExpl}), to the spaces ${\cal
S}^k_{\zd}(\R^{2n+1})$ of symmetric contravariant density valued
tensor fields over $\R^{2n+1}$. This generalized Hamiltonian
\[X: {\mathcal S}^k_\delta(\R^{2n+1})\raa{\mathcal
S}^{k+1}_{\delta+\op{I}}(\R^{2n+1})\] maps $S=S(\xi)$ to \be
X(S)=X(S)(\xi)=\sum_j(\xi_{q^j}\partial_{p_j}-
\xi_{p_j}\partial_{q^j})S(\xi)+\xi_t{\cal E}_sS(\xi)-\langle {\cal
E}_s,\xi_s\rangle
\partial_tS(\xi)+a(k,\delta)\xi_tS(\zx),\label{GenHam}\ee where
$a(k,\delta)=2(n+1)\delta-k.$ If $k=0$, operator $X$ obviously
coincides with the map $$X:{\cal S}^0_{\zd}(M)\ni
f\to\zs_1(\{f,.\})\in{\cal S}^1_{\zd+\op{I}}(M),$$ where $\{.,.\}$
is the Lagrange bracket and where $\{f,.\}:{\cal S}^0_0(M)\ni
g\to\{f,g\}\in{\cal S}^{0}_{\zd+\op{I}}(M)$. In \cite{FMP}, we
proved the following
\begin{prop}
Operator $X :{\mathcal S}^k_\delta(\R^{2n+1})\to{\mathcal
S}^{k+1}_{\delta+\op{I}}(\R^{2n+1})$ intertwines the
$\op{sp}_{2n+2}$-action and does not commute with the
$\op{CVect}(\R^{2n+1})$-action, unless $k=0$.
\label{XIntertwAct}\end{prop} Since the operators $i_{\alpha}$ and
$X$ modify the density weight of their arguments, we introduce the
space
\[R_{\delta}=\oplus_{k\in\N}R^k_\delta:=\oplus_{k\in\N}{\cal
S}^k_{\delta+\frac{k}{n+1}}(\R^{2n+1}).\]

\begin{theo}\label{ClassSpInvSame} The algebra of $\op{sp}_{2n+2}$-invariant
$($ differential $)$ operators from $R_{\delta}$ into $R_{\delta}$
is generated by $i_{\alpha}$ and $X$. More precisely, the space of
$\op{sp}_{2n+2}$- invariant operators from $R^{\E}_{\delta}$ into
$R^k_{\delta}$ is spanned by $\{X^m\circ i_{\za}^{\E+m-k} :
\op{sup}(0,k-\E)\leq m\leq k\}.$
\end{theo}

{\it Proof}. In order to simplify notations, we set
$M:=\R^{2n+1}$. Let $I$ be an $\op{sp}_{2n+2}$-invariant
differential operator, say of order $m$, from $R^{\E}_{\delta}$
into $R^k_{\delta}$. The principal symbol $\zs_m(I)$ of $I$ is an
invariant tensor field in
\[\Gamma({\cal S}^mTM\otimes \op{Hom}({\cal S}^{\E}TM\otimes
\mathbb{F}_{\delta+\frac{\E}{n+1}}M,{\cal
S}^{k}TM\otimes\mathbb{F}_{\delta+\frac{k}{n+1}}M))\simeq\Gamma({\cal
S}^kTM\otimes {\cal S}^mTM\otimes {\cal S}^{\E}T^*M\otimes
\mathbb{F}_{\frac{k-\E}{n+1}}M),\] see Remark 2 below. In view of
Theorem \ref{sp} and Theorem \ref{invcont}, we then have
\[\sigma_m(I)=\sum_{a,b,c,d\in\N}C_{abcd}\;u_1^{a}u_2^bu_3^cL_1^d\quad\quad(C_{abcd}\in\R),\]
where $a,b,c,$ and $d$ are subject to the conditions
\[\left\{\begin{array}{lll}
a+d&=&k\\
b+d&=&m\\
a+b+c&=&\E\\
d-c&=&k-\E
\end{array}\right..\] This system has a unique solution $a=k-m$, $b=0$,
$c=\E+m-k$, and $d=m$, with $\op{sup}(0,k-\E)\leq m\leq k$. Hence,
\[\sigma_m(I)(\eta,\xi;Y^{\E})=C\;u_1^{k-m}u_3^{\E+m-k}L_1^m\quad\quad(C\in\R),\] where we used conventional notations of affine
symbol calculus, see Remark 2 below.

Observe now that operator $X^m\circ
i_{\alpha}^{\E+m-k}:R^\E_{\delta}\to R^k_{\delta}$ is an
$\op{sp}_{2n+2}$-invariant differential operator of order $m$
(since $X$ has order 1 and $i_{\alpha}$ has order 0). Thus its
principal symbol reads
\[\sigma_m(X^m\circ i_{\alpha}^{\E+m-k})(\eta,\xi;Y^\E)=C_0\; u_1^{k-m}u_3^{\E+m-k}L_1^m\quad\quad(C_0\in\R_0).\]

It follows that the operator $I-\frac{C}{C_0}X^m\circ
i_{\alpha}^{\E+m-k}:R^\E_{\zd}\to R^k_{\zd}$ is
$\op{sp}_{2n+2}$-invariant and of order $\le m-1$. An easy
induction on the order of differentiation then
yields the result. \rule{1.5mm}{2.5mm}\\

{\bf Remark 2}. a. Let us mention that operator $X^m$ is tightly
related with the $m$th order Lagrange bracket. This
observation will be further developed in a subsequent work.\\

b. Although it is a known result in equivariant quantization,
commutativity---for differential operators between symmetric
contravariant density valued tensor fields---of the principal
symbol and the canonical actions of vector fields might not be
obvious for all the readers. Beyond computations in local
coordinates, affine symbol calculus allows to elegantly make sure
of the validity of this statement. Affine symbol calculus is a
non-standard computing technique. For further information we refer
the interested reader to \cite{NP}. Below we give the proof, via
symbol calculus, of the aforementioned commutativity.

Let $T\in\op{Hom}({\cal S}^{k'}_{\zd'}(\R^p),{\cal
S}^{k''}_{\zd''}(\R^p))_{\op{loc}}$, where subscript
``$\op{loc}$'' means that we confine ourselves to support
preserving operators. We fix coordinates and call {\it affine
symbol} $\zs_{\op{aff}}(T)$ of $T$, its total symbol (the highest
order terms of which coincide with the principal symbol $\zs(T)$).
Hence, if $T$ has order $m$, $\zs_{\op{aff}}(T)\in\zG({\cal
S}_mT\R^p\otimes\op{Hom}({\cal
S}^{k'}T\R^p\otimes\mathbb{F}_{\zd'}\R^p,{\cal
S}^{k''}T\R^p\otimes\mathbb{F}_{\zd''}\R^p))$, where ${\cal
S}_mT\R^p$ is the $m$th order filter (of the increasing
filtration) associated with the natural grading of ${\cal
S}T\R^p$. It is easily checked that, for any
$X\in\op{Vect}(\R^p)$,
\be\begin{array}{lll}\zs_{\op{aff}}(L_XT)&=&(X.T)(\zh;Y^{k'})-\la
X,\zh\ra\zt_{\zy}\lp
T(\zh;Y^{k'})\rp-X(\zy\p_{\xi})T(\zh;Y^{k'})\\
&&+T(\zh+\zy;X(\zy\p_{\xi})Y^{k'})+\zd''\la X,\zy\ra
T(\zh;Y^{k'})-\zd'\la X,\zy\ra
T(\zh+\zy;Y^{k'}),\end{array}\label{affsymbLXclassinv}\ee where we
used standard notations, see \cite{NP} ($X.T$ denotes the
derivatives of the coefficients of $T$, $Y^{k'}=Y\vee\ldots\vee Y$
($k'$ factors), $\zt_{\zy}\star$ is just a notation for the
translation $\star(.+\zy)-\star(.)$, $\zh$ symbolizes the
derivatives that act on the argument---represented by
$Y^{k'}$---of $T$, $\zy$ symbolizes the derivatives of the
coefficients of $X$, and $\xi$ denotes the variable of the
polynomials $Y^{k'}$ and $T(\zh;Y^{k'})$)%, and where the first
% term of the second line can be written
% $T(\zh+\zy;X(\zy\p_{\xi})Y^{k'})=\la Y,\zy\ra
% (X\p_Y)T(\zh+\zy;Y^{k'})$
. On the other hand, we have
$$\begin{array}{ll} L_X\zs (T)=&(X.\zs (T))(\zh;Y^{k'})-\la
X,\zh\ra(\zy\p_{\zh})\zs (T)(\zh;Y^{k'})\\&-X(\zy\p_{\xi})(\zs
(T)(\zh;Y^{k'}))+\zs (T)(\zh;X(\zy\p_{\xi})Y^{k'})+(\zd''-\zd')\la
X,\zy\ra \zs (T)(\zh;Y^{k'}).\end{array}$$ When selecting the
highest order terms in Equation (\ref{affsymbLXclassinv}), we see
that $\zs (L_XT)=L_X\zs (T),\forall X\in\op{Vect}(\R^p).$

A similar proof is possible for differential operators $T$ acting
between tensor densities. The corresponding result has already
been used earlier in this note. However, the observation that the
principal symbol intertwines the actions by Lie derivatives on
operators and symbols, is not true in general. It is for instance
not valid for ``quantum level operators'' $T\in\op{Hom}({\cal
D}^{k'}_{\zl'\zm'}(\R^p),{\cal
D}^{k''}_{\zl''\zm''}(\R^p))_{\op{loc}}.$

\section{Casimir operator}\label{CasimirSection}

As an application of our decomposition of the module of symbols
into submodules, see \cite{FMP}, of Section
\ref{InfProjContTrans}, and Section \ref{InvOpSameDens}, we now
prove that the Casimir operator $C^k_{\zd}$ of the canonical
representation of $\op{sp}_{2n+2}$ on $R^k_{\zd}$ (with respect to
the Killing form) is diagonal. Computation of this Casimir is a
challenge by itself, but, in addition, it will turn out that this
operator imposes restrictions on the parameters $k,k',\zd,\zd'$ of
symbol modules $R^k_{\zd}$ and $R^{k'}_{\zd'}$ that are
implemented by different density weights $\zd,\zd'$ and are linked
by an invariant operator.\\

The following upshots are well-known and mostly easily checked.
The symplectic algebra $\op{sp}(2n,\C)$ is a classical simple Lie
algebra of type $C_n$ (if $n\geq 3$). Its Killing form $K$ reads
$K: \op{sp}(2n,\C)\times \op{sp}(2n,\C)\ni(S,S')\raa
2(n+1)\op{tr}(SS')\in\C$, and its classical Cartan subalgebra
$C\subset\op{sp}(2n,\C)$ is $C=\{\op{diag}(\zD,-\zD),
\zD=\op{diag}(\zD_1,\ldots,\zD_n),\zD_i\in\C\}.$ The corresponding
roots are \be -(\delta_i+\delta_j)\quad (i\le j\le n),\quad
\delta_i+\delta_j\quad (i\le j\le n),\quad
\mbox{and}\quad\delta_i-\delta_j\quad (i,j\leq n),\label{roots}\ee
where $\delta_k$ is the $\C$-linear form of $C$ defined by
$\delta_k(\op{diag}(\Delta,-\Delta))=\Delta_k.$ If
$(e_1,\ldots,e_{2n})$ denotes as above the canonical basis of
$\C^{2n}$ and $(\epsilon^1,\ldots,\epsilon^{2n})$ the dual basis
in $\C^{2n*}$, the respective eigenvectors are
\begin{equation}\label{eigenvectors}
\epsilon^j\otimes e_{i+n}+\epsilon^i\otimes
e_{j+n},\,-\epsilon^{j+n}\otimes e_i-\epsilon^{i+n}\otimes e_j,\,
\mbox{and}\,-\epsilon^j\otimes e_i+\epsilon^{i+n}\otimes
e_{j+n}.\end{equation} We thus recover the result that the
eigenspaces $\op{sp}_{\zd}$ associated with the above-detailed
roots $\zd$, see Equation (\ref{roots}), are $1$-dimensional for
$\zd\neq 0$. Moreover, if $\Lambda$ denotes the set of roots, we
have the decomposition
$\op{sp}(2n,\C)=\bigoplus_{\zd\not=0,\zd\in\Lambda}\op{sp}_{\zd}\,\oplus
\,\,C.$ This splitting allows computing the Killing-dual basis of
basis (\ref{eigenvectors}), see also Equation (\ref{basis}).

\begin{prop}\label{DualBasisMatr}
The bases
\[\epsilon^j\otimes e_{i+n}+\epsilon^i\otimes
e_{j+n}\,(i\le j\le n),\,-\epsilon^{j+n}\otimes
e_i-\epsilon^{i+n}\otimes e_j\,(i\le j\le n),\,-\epsilon^j\otimes
e_i+\epsilon^{i+n}\otimes e_{j+n}\,(i,j\le n)\] and
\[k_{ij}(-\epsilon^{j+n}\otimes
e_{i}-\epsilon^{i+n} \otimes e_{j}),\,k_{ij}(\epsilon^{j}\otimes
e_{i+n}+\epsilon^{i}\otimes e_{j+n}),\,k(-\epsilon^i\otimes
e_j+\epsilon^{j+n}\otimes e_{i+n})\] of $\op{sp}(2n,\C)$ are dual
with respect to the Killing form, if and only if
$k_{ij}=-1/(4(n+1)(1+\delta_{ij}))$ and $k=1/(4(n+1))$.
\end{prop}

{\it Proof}. Remember first that if $N$ is a nilpotent subalgebra
of a complex Lie algebra $L$, and if $\zl,\zm\in \zL$ are roots of
$N$, such that $\zm\neq -\zl$, then the corresponding eigenspaces
$L_{\zl}$ and $L_\zm$ are orthogonal with respect to the Killing
form $K$ of $L$. Further, the basis $-\epsilon^i\otimes
e_i+\epsilon^{i+n}\otimes e_{i+n}$ $(i\in\{1,\ldots,n\})$ of $C$
is orthogonal with respect to $K$. Hence, it suffices to compute
$K$ on each pair of nonorthogonal vectors. For instance, we have
\[\begin{array}{l}
k_{ij}K(-\epsilon^{j+n}\otimes e_i-\epsilon^{i+n}\otimes
e_j,\epsilon^{j}\otimes e_{i+n}+\epsilon^{i}\otimes e_{j+n})
\\=-2(n+1)k_{ij}\op{tr}((\epsilon^{j+n}\otimes
e_i+\epsilon^{i+n}\otimes e_j)(\epsilon^{j}\otimes
e_{i+n}+\epsilon^{i}\otimes
e_{j+n}))\\
=-2(n+1)k_{ij}\op{tr}(\delta_{ij}\epsilon^i\otimes e_i +
\epsilon^i\otimes e_i + \epsilon^j\otimes
e_j+\delta_{ij}\epsilon^i\otimes e_i)\\ =1.
\end{array}\] The result follows.  \rule{1.5mm}{2.5mm}\\

{\bf Remarks.} \begin{itemize}\item All the matrices used above
are actually real matrices. The result on Killing-dual bases still
holds true for $\op{sp}(2n,\R)$ ($\op{sp}(2n,\R)$ is a split real
form of $\op{sp}(2n,\C)$, the Killing form of $\op{sp}(2n,\R)$ is
the restriction of the Killing form of $\op{sp}(2n,\C)$).\item If
read through Lie algebra isomorphism $J^{-1}\circ{\cal
X}:(\op{Pol}^2(\R^{2n}),\{.,.\}_{\zP})\raa
(\op{sp}(2n,\R),[.,.]_{\circ})$, see Equation (\ref{BasisPoly}),
Proposition \ref{DualBasisMatr} states that the bases
$$p_ip_j\;(i\le j\le n),\,q^{i}q^j\;(i\le j\le n),\,
p_jq^{i}\;(i,j\le n)$$ and $$k_{ij}\,q^{i}q^j\;(i\le j\le
n),\,k_{ij}\,p_ip_j\;(i\le j\le n),\, k\,p_iq^{j}\;(i,j\le n)$$ of
$\op{Pol}^2(\R^{2n})$ are Killing-dual. \item The preceding
result, written for space $\R^{2n+2}$ (coordinates:
$(p_1,\ldots,p_n,q^1,\ldots,q^n;t,\zt)$) and read through Lie
algebra isomorphism $X\circ\chi:
(\op{Pol}^2(\R^{2n+2}_+),\{.,.\}_{\zP})\raa
(\op{sp}_{2n+2},[.,.])$, see Equations (\ref{Chi}) and
(\ref{ContHamMapExpl}), shows that the bases \be\begin{array}{c}
X_{p_ip_j}\;(i\le j\le
n),\,X_{tp_i}\;(i\in\{1,\ldots,n\}),\,X_{t^2};\\X_{q^{i}q^j}\;(i\le
j\le n),\, X_{q^{i}}\;(i\in\{1,\ldots,n\}),\,X_1;\\
X_{p_jq^{i}}\;(i,j\le
n),\,X_{tq^{i}}\;(i\in\{1,\ldots,n\}),\,X_{p_i}\;(i\in\{1,\ldots,n\}),\,X_t
\end{array}\label{BasisSpFin}\ee and
\be\begin{array}{c}k_{ij}\,X_{q^{i}q^j}\;(i\le j\le
n),\,-k\,X_{q^{i}}\;(i\in\{1,\ldots,n\}),\,
-k/2\,X_1;\\k_{ij}\,X_{p_ip_j}\;(i\le j\le
n),\,-k\,X_{tp_i}\;(i\in\{1,\ldots,n\}),\,-k/2\,X_{t^2};\\k\,X_{p_iq^{j}}\;(i,j\le
n),\,k\,X_{p_i}\;(i\in\{1,\ldots,n\}),\,k\,X_{tq^{i}}\;(i\in\{1,\ldots,n\}),\,k\,X_t,\end{array}\label{BasisSpDual}\ee
with $k_{ij}=-1/(4(n+2)(1+\delta_{ij}))$ and $k=1/(4(n+2))$, are
bases of the algebra $\op{sp}_{2n+2}$ of infinitesimal projective
contact transformations, which are dual with respect to the
Killing form. Observe that the first basis is the basis computed
in Section \ref{InfProjContTrans} and that both bases are
explicitly known, see Equations (\ref{BasisSp1}),
(\ref{BasisSp2}), (\ref{BasisSp3}), and (\ref{BasisSp4}).
\end{itemize}

We already mentioned that action (\ref{ActionSymb}) of
$X\in\op{Vect}(\R^{2n+1})$ on $P\in {\cal S}_{\zd}(\R^{2n+1})$ has
the explicit form \be
L_XP=X(P)-\p_jX^{i}\xi_i\p_{\xi_j}P+\zd\,\p_iX^{i}\,P,\label{ActionSymbExpl}\ee
see Equation (\ref{invcond}). Remark that in this section we
denote the base coordinates by $(p_1,\ldots,p_n,q^1,\ldots,q^n,t)$
and the fiber coordinates by
$(\xi_{p_1},\ldots,\xi_{p_n},\xi_{q^1},\ldots,\xi_{q^n},\xi_{t})$.
Moreover, we took an interest in the Casimir operator $C^k_{\zd}$
of the preceding action of $\op{sp}_{2n+2}$ on $R^k_{\zd}={\cal
S}^k_{\zd+k/(n+1)}(\R^{2n+1})$, so that the weight in Equation
(\ref{ActionSymbExpl}) must be modified accordingly. The actions
on $R^k_{\zd}$ of the dual bases (\ref{BasisSpFin}) and
(\ref{BasisSpDual}) are now straightforwardly obtained:
\begin{equation}\label{der}\begin{array}{lll}
L_{X_{1}}&=& - 2\p_{t},\\
L_{X_{p_i}}&=&\p_{q^i} - p_i\p_{t} + \xi_t\p_{\xi_{p_i}},\\
L_{X_{q^i}}&=&-\p_{p_i} - q^i\p_{t}+ \xi_t\p_{\xi_{q^i}},\\
L_{X_{t}}&=& - {\cal E}_s- 2t\p_{t}+{\cal E}_{\xi_s} +
2\xi_t\p_{\xi_t} - 2((n+1)\zd+k),\\
L_{X_{p_ip_j}}&=&p_j\p_{q^{i}}+p_i\p_{q^j}-\xi_{q^{i}}\p_{\xi_{p_j}}-{\xi_{q^j}}\p_{\xi_{p_i}},\\
L_{X_{q^{i}q^{j}}}&=&-q^j\p_{p_i}-q^{i}\p_{p_j}+\xi_{p_i}\p_{\xi_{q^j}}+\xi_{p_j}\p_{\xi_{q^{i}}},\\
L_{X_{p_jq^{i}}}&=&q^{i}\p_{q^j}-p_j\p_{p_i}+\xi_{p_i}\p_{\xi_{p_j}}-\xi_{q^j}\p_{\xi_{q^{i}}},\\
L_{X_{tp_i}}&=&t(\p_{q^i} - p_i\p_{t})- p_i{\cal E}_s -
\xi_{q^i}\p_{\xi_t}+p_i{\cal E}_{\xi} + {\cal
E}(\xi)\p_{\xi_{p_i}}-2((n+1)\zd+k) p_i,
\\
L_{X_{tq^i}}&=&-t(\p_{p_i}+q^i\p_{t})-q^i{\cal E}_s+
\xi_{p_i}\p_{\xi_t}+q^i{\cal E}_{\xi}
 + {\cal E}(\xi)\p_{\xi_{q^i}}-2((n+1)\zd+k) q^i,\\
L_{X_{t^2}}&=&-2t{\cal E}+2t{\cal E}_{\xi}+2{\cal
E}(\xi)\p_{\xi_t}-4((n+1)\zd+k) t.
\end{array}
\end{equation}
In these equations, ${\cal E}$ is the Euler field of $\R^{2n+1}$,
${\cal E}_s$ is its spatial part, ${\cal E}_{\xi}$ is the Euler
field with respect to the fiber coordinates, ${\cal
E}_{\xi}=\xi_{p_i}\p_{\xi_{p_i}}+\xi_{q^{i}}\p_{\xi_{q^{i}}}+\xi_t\p_{\xi_t},$
${\cal E}_{\xi_s}$ denotes the spatial part of ${\cal E}_{\xi}$,
and ${\cal E}(\xi)$ is the contraction of ${\cal E}$ and
$\xi=\xi_{p_i}dp_i+\xi_{q^{i}}dq^{i}+\xi_tdt$.\\

When combining Equations (\ref{BasisSpFin}), (\ref{BasisSpDual}),
and (\ref{der}), we get the explicit form of Casimir operator
$C^k_{\zd}$. A direct computation and simplification of
$C^k_{\zd}$ are possible, but almost inextricable. In the sequel,
we provide a much more economic method based upon
Theorem \ref{ClassSpInvSame}.\\

Set ${\frak C}_k=\{-p/(2(n+1)):p=0,1,\dots,2k-2\}$.

\begin{theo}\label{CasDiagForm} If $\zd\notin {\frak C}_k$, the Casimir operator $C^k_{\zd}$ of the Lie algebra $\op{sp_{2n+2}}$ of infinitesimal
projective contact transformations, with respect to its Killing
form, and for its canonical action (\ref{ActionSymbExpl}) on
$R^k_\delta={\cal S}^k_{\zd+k/(n+1)}(\R^{2n+1})$, is given by
\[C^k_{\zd}=\frac{1}{n+2}(c(k,\delta)\op{id}+X\circ i_\alpha),\]
where
\[c(k,\delta)= (n+1)^2\delta^2-(n+1)^2\delta+k(n+1)\delta+\frac{k^2-k}{2}\]
and where $X$ $($resp. $i_{\za}$$)$ is the generalized Hamiltonian
$($resp. the vertical cotangent lift of $\za$$)$ defined in
Section \ref{InvOpSameDens}.
\end{theo}

Let us first recall two results obtained in \cite{FMP}.

\begin{prop}\label{ComzaX^m} In $R^k_{\zd}$ and for $\E\in\N_0,$ we have $$i_\za\circ X^{\E} = X^{\E}\circ
i_\za+r(\E,k)X^{\E-1},$$ where
$r(\E,k)=-\frac{\E}{2}(2(n+1)\zd+2k+\E-1)$.\end{prop}

\begin{theo}\label{DecSubMod} If $\zd\notin{\frak C}_k,$ then $$R^k_{\zd}=\bigoplus_{m=0}^kR_{\zd}^{k,\E}:=
\bigoplus_{\E=0}^kX^{\E}(R^{k-\E}_{\zd}\cap\op{ker}i_{\za}).$$\end{theo}

The last upshot, which extends splitting (\ref{OvsDecomp}), is the
main result of \cite{FMP} and has actually been proved in view of
the present application.\\

{\it Proof of Theorem \ref{CasDiagForm}}. The proof consists of three stages.\\

1. Casimir operator $C^k_{\zd}$ reads \be\begin{array}{lll}
C^k_{\zd}&=&-\frac{1}{8(n+2)}(L_{X_{1}}\circ L_{X_{t^2}}
+L_{X_{t^2}}\circ L_{X_{1}}) +\frac{1}{4(n+2)}(L_{X_{t}})^2
-\frac{1}{4(n+2)}\sum_i(L_{X_{q^i}}\circ
L_{X_{tp_i}}+L_{X_{tp_i}}\circ L_{X_{q^i}})\\
&+&\frac{1}{4(n+2)}\sum_i(L_{X_{p_i}}\circ L_{X_{tq^i}}+
L_{X_{tq^i}}\circ L_{X_{p_i}})
-\frac{1}{8(n+2)}\sum_i(L_{X_{q^{i2}}}\circ L_{X_{p_i^2}} +L_{X_{p_i^2}}\circ L_{X_{q^{i2}}})\\
&-&\frac{1}{4(n+2)}\sum_{i<j}(L_{X_{q^i q^j}}\circ L_{X_{p_i
p_j}}+L_{X_{p_i p_j}}\circ L_{X_{q^i q^j}})
+\frac{1}{4(n+2)}\sum_i\sum_j L_{X_{p_iq^j}}\circ L_{X_{p_jq^i}}.
\end{array}\label{Casimir1}\ee
Since
$X\circ\chi:(\op{Pol}^2(\R^{2n+2}_+),\{.,.\}_{\zP})\raa(\op{sp}_{2n+2},[.,.])$
and
$L:(\op{sp_{2n+2},[.,.]})\raa(\op{End}(R^k_{\zd}),[.,.]_{\circ})$
are Lie algebra homomorphisms, and $\{\tau^2,t^2\}_{\zP}=-4t\tau$,
we have $[L_{X_{1}},L_{X_{t^2}}]_{\circ}=-4L_{X_t}$. Hence, the
first term of the {\small RHS} of Equation (\ref{Casimir1}) is
equal to $-1/(4(n+2))(L_{X_{t^2}}\circ L_{X_1}-2L_{X_t})$. When
using similarly the Poisson brackets $\{ \tau q^i,t p_i
\}_{\zP}=-p_iq^i-t\tau,$ $\{ \tau p_i,t q^i \}_{\zP}=-p_iq^i+t
\tau,$ $\{q^{i2},p^2_{i}\}_{\zP}=-4p_iq^{i}$,
$\{q^{i}q^j,p_ip_j\}_{\zP}=-p_iq^{i}-p_jq^j$ ($i\neq j$), we
finally get \be\begin{array}{lll} C^k_{\zd}&=&-\frac{1}{4(n+2)}
L_{X_{t^2}}\circ L_{X_{1}}+\frac{1}{4(n+2)}(L_{X_{t}})^2
+\frac{1}{2(n+2)}\sum_i(L_{X_{tq^i}}\circ
L_{X_{p_i}}-L_{X_{tp_i}}\circ L_{X_{q^i}})\\
&-&\frac{1}{4(n+2)}\sum_iL_{X_{p_i^2}}\circ L_{X_{q^{i2}}}
-\frac{1}{2(n+2)}\sum_{i<j}L_{X_{p_ip_j}}\circ L_{X_{q^{i}q^j}}\\
&+&\frac{1}{4(n+2)}\sum_i\sum_jL_{X_{p_iq^j}}\circ
L_{X_{p_jq^{i}}}+\frac{1}{2}\frac{n+1}{n+2}L_{X_{t}}+\frac{n+1}{4(n+2)}\sum_i
L_{X_{p_iq^{i}}}.
\end{array}\label{Casimir2}\ee

2. Theorem \ref{ClassSpInvSame} entails that \be
C^k_{\zd}=\sum_{m=0}^{\op{inf}(2,k)}c^k_{\zd,m}\;X^m\circ
i_{\za}^m\quad (c^k_{\zd,m}\in\R),\label{CasInv}\ee where we have
used the fact that Casimir $C^k_{\zd}$ is a second order
differential operator, see Equation (\ref{ActionSymbExpl}).
Observe also that it follows from Proposition \ref{ComzaX^m} that
$R^{k,\E}_{\zd}=X^{\E}(R^{k-\E}_{\zd}\cap\op{ker}i_{\za})$,
$\E\in\{0,\ldots,k\}$, is an eigenspace of $C_{\zd}^k$ with
eigenvalue
\be\ze^{k,\E}_{\zd}=\sum_{m=0}^{\op{inf}(2,\E)}c^k_{\zd,m}\zP_{i=\E-m+1}^{\E}r(i,k-\E).\label{eigenvalues}\ee
In particular, \be\begin{array}{c}
C^k_{\zd}\vert_{R^{k,0}_{\zd}}=c^k_{\zd,0}\op{id},C^k_{\zd}\vert_{R^{k,1}_{\zd}}=(c^k_{\zd,0}+r(1,k-1)c^k_{\zd,1})\op{id},\\
C^k_{\zd}\vert_{R^{k,2}_{\zd}}=(c^k_{\zd,0}+r(2,k-2)c^{k}_{\zd,1}+r(1,k-2)r(2,k-2)c^k_{\zd,2})\op{id}.\label{CasEigen}\end{array}\ee
%When computing the Casimir operator successively on elements of
%$R^{k\E}_{\zd}=X^{\E}(R^{k-\E}_{\zd}\cap\op{ker}\za)$
%($\E\in\{0,1,2\}$), using as well Equation (\ref{CasInv}) as
%Equation (\ref{CasExpl}), we get a triangular system in the
%unknowns $c^k_{\zd,m}$ ($m\in\{0,1,2\}$).
%Observe that the projectors allow to obtain elements in
%$R^{k\E}_{\zd}=X^{\E}(R^{k-\E}_{\zd}\cap\op{ker}\za)$, but that
%such elements can also be obtained by hand. For instance,
%$(p_1\xi_t+\xi_{q^1})^k\in R^{k0}_{\zd}$ (other elements are
%$(q^1\xi_t-\xi_{p_1})^k$, $(q^1\xi_{q^1}+p_1\xi_{p_1})^k$),
%$X(p_1\xi_t+\xi_{q^1})^{k-1}\in R^{k1}_{\zd}$, ...\\
%Now, in principle, the Casimir is computed (maybe some problems
%appear for vanishing coefficients).

Let now $P^k_{\zd_a}$ be the polynomial $(p_1\xi_t+\xi_{q^1})^k$
viewed as an element of ${\cal
S}^k_{\zd+\frac{a}{n+1}}(\R^{2n+1})$. Since, for any $S\in
R^k_{\zd}$, the contraction $i_{\za}S\in R^{k-1}_{\zd}$ is given
by
$$(i_{\za}S)(\xi)=\frac{1}{2}\left(\sum_j(p_j\p_{\xi_{q^j}}-q^j\p_{\xi_{p_j}})-\p_{\xi_t}\right)S(\xi),$$ see Section
\ref{InvOpSameDens}, it is clear that $P^k_{\zd_k}\in
R^{k,0}_{\zd}.$ Moreover, it is easily checked, see Equation
(\ref{GenHam}), that
$$X(P^{k-1}_{\zd_{k-1}})=(2(n+1)\zd+2k-2))\xi_tP^{k-1}_{\zd_{k}}=-2r(1,k-1)\xi_tP^{k-1}_{\zd_{k}}\in
R^{k,1}_{\zd}$$ and that
$$X^2(P^{k-2}_{\zd_{k-2}})=(2(n+1)\zd+2k-4)(2(n+1)\zd+2k-3)\xi_t^2P^{k-2}_{\zd_{k}}=2r(1,k-2)r(2,k-2)\xi_t^2P^{k-2}_{\zd_{k}}\in
R^{k,2}_{\zd},$$ where the coefficients in the {\small RHS} of the
two last equations do not vanish, since $\zd\notin{\frak C}_k$. As
generalized Hamiltonian $X$ intertwines the
$\op{sp}_{2n+2}\,$-action, see Proposition \ref{XIntertwAct}, we
also have
$C^k_{\zd}X(P^{k-1}_{\zd_{k-1}})=XC^{k-1}_{\zd}(P^{k-1}_{\zd_{k-1}})$.
If we now apply Equation (\ref{CasEigen}) to both sides, we get
$(c^k_{\zd,0}+r(1,k-1)c^k_{\zd,1})X(P^{k-1}_{\zd_{k-1}})=c^{k-1}_{\zd,0}X(P^{k-1}_{\zd_{k-1}}).$
Thus, \be
c^k_{\zd,1}=\frac{c^{k-1}_{\zd,0}-c^k_{\zd,0}}{r(1,k-1)}.\label{C1}\ee
When proceeding analogously for $X^2(P^{k-2}_{\zd_{k-2}}),$ we
obtain \be
c^{k}_{\zd,2}=\frac{r(1,k-1)(c^{k-2}_{\zd,0}-c^{k}_{\zd,0})-r(2,k-2)(c^{k-1}_{\zd,0}-c^k_{\zd,0})}{r(1,k-1)r(1,k-2)r(2,k-2)}.\label{C2}\ee

3. Hence, Casimir operator $C^k_{\zd}$ is completely known, see
Equations (\ref{CasInv}), (\ref{C1}), and (\ref{C2}), if we find
$c^k_{\zd,0}$. In this effect, we use Equation (\ref{CasEigen})
for $P^k_{\zd_k}\in R^{k,0}_{\zd}$, and compute the {\small LHS}
by means of Equation (\ref{Casimir2}). Straightforward (and even
fairly short) computations allow checking the contributions of the
successive terms $T_1$--$T_8$ of the {\small RHS} of
(\ref{Casimir2}).
$$\begin{array}{c}T_1P^k_{\zd_k}=0,T_2P^k_{\zd_k}=\frac{1}{4(n+2)}(2(n+1)\zd+k)^2P^k_{\zd_k},
T_3P^k_{\zd_k}=0, T_4P^k_{\zd_k}=\frac{k}{n+2}P^k_{\zd_k},
T_5P^k_{\zd_k}=\frac{k(n-1)}{2(n+2)}P^k_{\zd_k},\\
T_6P^k_{\zd_k}=\frac{k(n-1+k)}{4(n+2)}P^k_{\zd_k},T_7P^k_{\zd_k}=-\frac{1}{2}\frac{n+1}{n+2}(2(n+1)\zd+k)P^k_{\zd_k},
T_8P^k_{\zd_k}=-\frac{k(n+1)}{4(n+2)}P^k_{\zd_k}.\end{array}$$

When summing up these terms, we get \be c^k_{\zd,0}=
\frac{1}{n+2}((n+1)^2\delta^2-(n+1)^2\delta+k(n+1)\delta+(k^2-k)/2),\label{C0}\ee
and when substituting in Equations (\ref{C1}) and (\ref{C2}), we
obtain \be
c^k_{\zd,1}=\frac{1}{n+2}\;\mbox{and}\;c^k_{\zd,2}=0.\label{C20}\quad\rule{1.5mm}{2.5mm}\ee

\begin{prop}\label{EigenvaluesDiff} If $\zd\notin{\frak C}_k$, space $R^{k}_{\zd}$ is the
direct sum of the eigenspaces $R^{k,\E}_{\zd}$,
$\E\in\{0,\ldots,k\}$, of Casimir operator $C^k_{\zd}$. The
corresponding eigenvalues are
$\ze^{k,\E}_{\zd}=1/(n+2)(c(k,\zd)+r(\E,k-\E))$, see Theorem
\ref{CasDiagForm} and Proposition \ref{ComzaX^m}, and eigenvalues
$\ze^{k,\E}_{\zd}$ associated with different $\E$ cannot
coincide.\end{prop}

{\it Proof}. The first assertion and the values of the
$\ze^{k,\E}_{\zd}$, $\E\in\{0,\ldots,k\}$, are direct consequences
of Theorem \ref{DecSubMod} and Equations (\ref{eigenvalues}),
(\ref{C0}), and (\ref{C20}). Assume now that, for $\E_1\neq\E_2$,
we have $\ze^{k,\E_1}_{\zd}=\ze^{k,\E_2}_{\zd}$. This means that
$r(\E_1,k-\E_1)=r(\E_2,k-\E_2)$, i.e. that
$(\E_1-\E_2)(2(n+1)\zd+2k-(\E_1+\E_2+1))=0$. As
$2\le\E_1+\E_2+1\le 2k$, the last result is possible only if
$\zd\in{\frak C}_k$.  \rule{1.5mm}{2.5mm}

\section{Invariant operators between symbol modules implemented by different density
weights}\label{InvOpDiffDens}

In this section, we use the above Casimir operator to show that
the quest for invariant operators $T:R^{k}_{\zd}\raa
R^{k'}_{\zd'}$ between symbol spaces (implemented by different
weights $\zd$ and $\zd'$) can be put down to the search of a small
number $\zk$ of invariant operators $T_{\E_i}:R^{k,\E_i}_{\zd}\raa
R^{k',\E'_i}_{\zd'}$, $i\in\{1,\ldots,\zk\}$, between (smaller)
eigenspaces (that are of course submodules the corresponding
modules $R^{k}_{\zd}$ and $R^{k'}_{\zd'}$).

\begin{prop} Let $T:R^{k}_{\zd}\raa R^{k'}_{\zd'}$ be an $\op{sp}_{2n+2}$-invariant
operator. Assume that $\zd\notin{\frak C}_k,$ $\zd'\notin{\frak
C}_{k'}$. For any eigenspace $R^{k,\E}_{\zd}$,
$\E\in\{0,\ldots,k\}$, the restriction $T_{\E}$ of $T$ to
$R^{k,\E}_{\zd}$, either vanishes, or is an
$\op{sp}_{2n+2}$-invariant operator from $R^{k,\E}_{\zd}$ into an
eigenspace $R^{k',\E'}_{\zd'}$, $\E'\in\{0,\ldots,k'\}$, where
$\E'$ verifies the equation
$\ze^{k',\E'}_{\zd'}=\ze^{k,\E}_{\zd}$. Moreover, the map
$':\E\raa \E'$ is injective.
\end{prop}

{\it Proof}. Since $T$ is invariant, we have
$TC^k_{\zd}=C^{k'}_{\zd'}T$, so that for any eigenvector $P\in
R^{k,\E}_{\zd}$, one gets
$$C^{k'}_{\zd'}TP=TC^k_{\zd}P=\ze^{k,\E}_{\zd}TP=\sum_{\E'=0}^{k'}\ze^{k,\E}_{\zd}p^{k',\E'}_{\zd'}TP,$$
where we denote by $p^{k',\E'}_{\zd'}:R^{k'}_{\zd'}\raa
R^{k',\E'}_{\zd'}$ the canonical projection onto submodule
$R^{k',\E'}_{\zd'}$. On the other hand,
$$C^{k'}_{\zd'}TP=\sum_{\E'=0}^{k'}C^{k'}_{\zd'}p^{k',\E'}_{\zd'}TP=\sum_{\E'=0}^{k'}\ze^{k',\E'}_{\zd'}p^{k',\E'}_{\zd'}TP.$$
If $T_{\E}$ does not vanish, there is at leat one $\E'$, such that
$p^{k',\E'}_{\zd'}T_{\E}\neq 0,$ and for this $\E'$, we have
$\ze^{k\E}_{\zd}=\ze^{k'\E'}_{\zd'}.$ Actually, $\E'$ is unique,
in view of Proposition \ref{EigenvaluesDiff}. Hence,
$T_{\E}:R^{k,\E}_{\zd}\raa R^{k',\E'}_{\zd'}$ (and is of course
invariant). Eventually, map $':\E\raa \E'$ is injective, again in
view of Proposition \ref{EigenvaluesDiff}.  \rule{1.5mm}{2.5mm}\\

In the following $T$ is again a nonvanishing
$\op{sp}_{2n+2}$-invariant operator between symbol spaces
$R^{k}_{\zd}$ and $R^{k'}_{\zd'}$, $\zd\notin{\frak C}_k,$
$\zd'\notin{\frak C}_{k'}$.

Let $T_{\E_i}:R^{k,\E_i}_{\zd}\raa R^{k',\E_i'}_{\zd'}$,
$\E_i\in\{0,\ldots,k\}$, $\E_i'\in\{0,\ldots,k'\}$,
$i\in\{1,\ldots,\zk\}$, be the nonvanishing restrictions of $T$.
We denote the condition $\ze^{k',\E_i'}_{\zd'}=\ze^{k,\E_i}_{\zd}$
by ${\cal R}_i$ and systematically substitute the system ${\cal
R}_1,{\cal R}_2'={\cal R}_1-{\cal R}_2,\ldots,{\cal
R}_{\zk}'={\cal R}_{1}-{\cal R}_{\zk}$ to the system ${\cal
R}_{1},\ldots,{\cal R}_{\zk}$. Indeed, a short computation shows
that Equation ${\cal R}_i$ reads
\be\begin{array}{l} 2(n+1)^2(\zd+\zd'-1)(\zd-\zd')+(k+k'-1)(k-k')+2(n+1)(k\zd-k'\zd')\\
-2(n+1)(\E_i\zd-\E'_i\zd')-2(\E_i
k-\E'_ik')+(\E_i+\E'_i+1)(\E_i-\E'_i)=0\quad(i\in\{1,\ldots,\zk\}).\end{array}\label{R}\ee
Hence, ${\cal R}_i$ is a polynomial equation in $k,k',\zd,\zd'$ of
(total) degree $2$ (with integer coefficients), whereas Equation
${\cal R}_j'$, i.e. \be 2(n+1)\zD_j\zd-2(n+1)\zD'_j\zd'+2\zD_j
k-2\zD'_jk'=\zD_j-\zD'_j+\zD_j\zS_j-\zD'_j\zS'_j=:\zl_j\quad(j\in\{2,\ldots,\zk\}),\label{R'}\ee
where $\zD_j=\E_1-\E_j,\zD'_j=\E'_1-\E'_j$ and
$\zS_j=\E_1+\E_j,\zS'_j=\E_1'+\E_j'$, is linear in $k,k',\zd,\zd'$
(with integer coefficients).

Observe that Equations (\ref{R}) and (\ref{R'}) are
Diophantine-type equations. Let us recall that Diophantine
equations are indeterminate polynomial equations with integer
variables. In other words, Diophantine equation systems define
algebraic surfaces and ask for lattice points on them. Y.
Matiyasevich's solution of Hilbert's 10th problem shows that there
is no algorithm that allows solving arbitrary Diophantine systems.
One of the most celebrated results in this field is of course A.
Wiles' conclusion concerning Fermat's Diophantine equation
$x^n+y^n=z^n$, $n>2$.

Equations (\ref{R}) and (\ref{R'}) contain the real variables
$\zd,\zd'$, as well as the integer variables $k,k'$, which---in
addition---have lower bounds.

First, it is clear that, if $k,k',\zd,\zd'$ are known, Equation
(\ref{R}), which is quadratic in $\E'$, allows computing for each
$\E\in\{0,\ldots,k\}$, the corresponding $\E'\in\{0,\ldots,k'\}$.
If no $\E$ generates an appropriate $\E'$, no (nontrivial)
invariant operator exists between the considered spaces.

But let us revert to our initial viewpoint and translate existence
of nonvanishing restrictions $T_{\E_i}:R^{k,\E_i}_{\zd}\raa
R^{k',\E'_i}_{\zd'}$ into {\it conditions} on $k,k',\zd,\zd'$. It
follows of course from the above remarks that our system of
equations cannot be completely solved in the considered general
setting. Nevertheless, there is evidence that possible invariant
operators between symbol spaces induced by different density
weights should not have many nonvanishing building blocks
$T_{\E_i}$.

The assumed existence of a nonvanishing invariant operator $T$
between the considered symbol spaces, entails that Equation
(\ref{R}) has at least one solution
$(\E_1,\E'_1)\in\{0,\ldots,k\}\times\{0,\ldots,k'\}$. Hence,
$\zd'$ is a {\it solution} of the corresponding quadratic equation
${\cal R}_1$, so that the discriminant $D$ of this equation is
positive. Since $D$ is quadratic in $\zd$ and as $\zd^2$ has
positive coefficient $16(n+1)^4$, weight $\zd$ is {\sl not
located} between the possible roots of $D=0$.

If operator $T$ has at least two nonvanishing building blocks
$T_{\E_1}$ and $T_{\E_2}$, i.e. if $\zk=2$, Equations ${\cal R}_1$
and ${\cal R}_2'$ hold true. The second equation provides $\zd'$
as {\it first order polynomial} in $k,k'$, and $\zd$. When
substituting in the first equation, we find that $\zd$ verifies an
at most quadratic equation $a_2\zd^2+a_1\zd+a_0=0$,
$a_i=a_i(k,k',\E_1,\E_2,\E'_1,\E'_2)$. It can be seen that this
equation {\it allows computing} $\zd$ or $k'$, and gives {\sl
further conditions} on the parameters.

Although explicit results are far too complicated to be written
down here, it should now be clear that $\zk>4$ entails a
contradiction.

It is even quite easy to understand that, in most cases, this
conclusion already holds for $\zk=4.$ Indeed, consider the
conditions ${\cal R}'_2$, ${\cal R}'_3$, and ${\cal R}'_4$. The
lines $l_j:=(2(n+1)\zD_j,-2(n+1)\zD'_j,2\zD_j,-2\zD'_j)$,
$j\in\{2,3,4\}$, of the coefficient-matrix of the {\small LHS} of
these equations are obviously linearly dependent. Hence, one line
is a linear combination of the two others, e.g.
$l_4=c_2l_2+c_3l_3,$ $c_2,c_3\in\R$, i.e.
$\zD_4=c_2\zD_2+c_3\zD_3$ and $\zD'_4=c_2\zD'_2+c_3\zD'_3$. It
follows from the structure of the {\small RHS} of ${\cal R}'_j$,
see Equation (\ref{R'}), that existence of the same linear
combination $\zl_4=c_2\zl_2+c_3\zl_3$ between the {\small RHS} of
${\cal R}'_2$, ${\cal R}'_3$, and ${\cal R}'_4$, is a strong
requirement.
% Actually, if we imagine the requirement as rank of coefficient-matrix = rank of coefficient matrix + RHS, we get a condition on 4 \E and the
% corresponding \E'=\E'(\E,k,k',\zd,\zd'), so a condition on $k,k',\zd,\zd'$ (I tend to think that as we substitute $R$ a condition
% that does not come from conditions $R'$, we get a nontrivial new condition, in opposition with what I observed in my explicit computations),
% i.e. a condition on $\zd'$.
If $\zl_4\neq c_2\zl_2+c_3\zl_3$, the considered system is
incompatible---a contradiction.

The last possible case, $\zk=3$, is also interesting. Indeed, if
$(\zD_j,\zD'_j)$, $j\in\{2,3\}$, are independent, and if we set
$$L=L(\zD_j,\zD_j',\zS_j,\zS_j')=\frac{\zD_2\zD'_3\zS_2-\zD_3\zD'_2\zS_3+\zD'_2\zD'_3(\zS'_3-\zS'_2)}{2(n+1)(\zD_2\zD'_3-\zD_3\zD'_2)}\in\Q,$$
it easily follows from the linear system ${\cal R}'_2$, ${\cal
R}'_3$ that
$$\zd=-\frac{2k-1}{2(n+1)}+L\in\Q.$$ A similar upshot $\zd'\in\Q$
is also valid.

\end{document}